\documentclass[12pt,draftcls,onecolumn]{IEEEtran}
\usepackage{bbm}
\usepackage{mathrsfs}
\usepackage{amsfonts}
\usepackage{epsfig}
\usepackage{psfrag}
\usepackage{graphicx}
\usepackage{epstopdf}
\usepackage{subfigure}
\usepackage{cite}

\usepackage{amssymb}
\usepackage{amsthm}
\usepackage{algorithm,algorithmic}
\usepackage[tbtags]{amsmath}
\usepackage{url}
\usepackage{bm}

\usepackage{color}
\definecolor{gray}{RGB}{128,128,128}

\allowdisplaybreaks

\newtheorem{theorem}{Theorem}
\newtheorem{assumption}{Assumption}

\newtheorem{lemma}{Lemma}

\newtheorem{remark}{Remark}

\newcommand{\lx}[1]{{\color{black} #1}}
\newcommand{\TY}[1]{{\color{black} #1}}
\newcommand{\red}[1]{{\color{black} #1}}
\newcommand{\blue}[1]{{\color{black} #1}}
\newcommand{\XL}[1]{{\color{black} #1}}
\newcommand{\xul}[1]{{\color{black} #1}}
\usepackage{color}
\definecolor{gray}{RGB}{128,128,128}

\allowdisplaybreaks

\begin{document}
\title{Compressed Gradient Tracking Algorithms for Distributed Nonconvex Optimization}
\author{Lei~Xu,~\IEEEmembership{Student Member, IEEE},
        ~Xinlei~Yi\IEEEmembership{},\\
        ~Guanghui~Wen,~\IEEEmembership{Senior Member, IEEE},
        ~Yang Shi,~\IEEEmembership{Fellow, IEEE},\\
        Karl~H.~Johansson,~\IEEEmembership{Fellow, IEEE},
        ~Tao~Yang,~\IEEEmembership{Senior Member, IEEE}
\thanks{This work was supported by the National Natural Science Foundation of China under Grants 62133003, 61991403 \& 61991400, the Knut and Alice Wallenberg Foundation, the
Swedish Foundation for Strategic Research.}
\thanks{L. Xu, and T. Yang are with the State Key Laboratory of Synthetical Automation for Process Industries, Northeastern University, Shenyang 110819, China(e-mail: 2010345@stu.neu.edu.cn; yangtao@mail.neu.edu.cn).(Corresponding author:
Tao~Yang.)}
\thanks{\red{X. Yi is with the Lab for Information \& Decision Systems, Massachusetts Institute of Technology, Cambridge, MA 02139, USA (e-mail: xinleiyi@mit.edu)}}
\thanks{\red{G. Wen is with the Department of Systems Science, School of Mathematics, Southeast University, Nanjing 210096, China (e-mail:
wenguanghui@gmail.com).}}
\thanks{\red{Y. Shi is with the Department of Mechanical Engineering, University of
Victoria, Victoria, BC V8W 2Y2, Canada (e-mail: yshi@uvic.ca).}}
\thanks{K. H. Johansson is with the Division of Decision and Control Systems, School of Electrical Engineering and Computer Science, KTH Royal Institute of Technology, 100 44, Stockholm, Sweden (e-mail: kallej@kth.se).}
}


\maketitle

\begin{abstract}
In this paper, we study the distributed nonconvex optimization problem, aiming to minimize the average value of the local nonconvex cost functions using local information exchange.
To reduce the communication overhead, we introduce three general classes of compressors, i.e., compressors with bounded relative compression error, compressors with globally bounded absolute compression error, and compressors with locally bounded absolute compression error.
By integrating them, respectively, with the distributed gradient tracking algorithm, we then propose three corresponding compressed distributed nonconvex optimization algorithms.
\XL{Motivated by the state-of-the-art BEER algorithm proposed in \cite{zhao2022NIPS}, which is an efficient compressed algorithm integrating gradient tracking with biased and contractive compressors, our first proposed algorithm extends this algorithm to accommodate both biased and non-contractive compressors.
}
\red{For each algorithm, we design a novel} Lyapunov function to demonstrate \red{its sublinear convergence} to a stationary point if the local cost functions are smooth.
Furthermore, when the global cost function satisfies the Polyak--{\L}ojasiewicz (P--{\L}) condition, we \red{show} that our proposed algorithms linearly converge to a global optimal point.
\red{It is} worth noting that, for compressors with bounded relative compression error and globally bounded \blue{absolute} compression error, our \red{proposed algorithms'} parameters do not require prior knowledge of the P--{\L} constant.
The theoretical results are illustrated by numerical example\TY{s}, demonstrating the effectiveness of the proposed algorithms in significantly reducing the communication burden while maintaining the convergence performance.
Moreover, \red{simulation results show} that the proposed algorithms outperform the state-of-the-art compressed distributed nonconvex optimization algorithms.
\end{abstract}


\IEEEpeerreviewmaketitle
\section{Introduction}
In recent years, distributed optimization has received considerable attention and has been applied in power systems, sensor networks, and machine learning, just to name a few \cite{Guo2016TPS,Du2018PESGM,Zhang2019TCNS,Sayed2014FTML}.
Various distributed optimization algorithms have been proposed; see, e.g., survey papers \cite{nedic2018distributed,Yang2019ARC} \red{and references therein}.
Early work \cite{Nedic2009TAC} proposed a distributed (sub)gradient descent (DGD) algorithm that requires a diminishing step size.
\red{In order to speed up the convergence rate}, the accelerated distributed optimization algorithms \red{with a fixed step size} have been proposed, such as EXTRA algorithm \cite{Shi2015SIAM}, distributed proportional-integral (DPI) algorithm \cite{wang2010AAC,Kia2015Aut}, and distributed gradient tracking (DGT) algorithm \cite{Nedic2017SIAM,Qu2017TCNS,Xu2015CDC}.

Distributed optimization algorithms involve data spread across multiple agents, which require each agent to communicate with its neighboring agents.
However, one of the main challenges is the communication bottleneck, which can occur owing to limited channel bandwidth or communication power.
Communication compression, which encompasses techniques like quantization and sparsification, can reduce the required communication capacity, see recent survey papers \cite{Shi2020CST,Cao2023JSAC}.
For distributed convex optimization problems,
various compressed distributed optimization algorithms have been developed.
\red{For example, based} on the DGD algorithm, \cite{Yi2014TCNS,Doan2020TAC,Zhang2018TAC} developed a quantized gradient algorithm using a uniform quantizer, a random quantizer, and the sign of the relative state\TY{, respectively}.
\cite{Alistarh2017NIPS,Horvath2022OMS} used a distributed stochastic gradient descent (DSGD) algorithm with gradient quantization and encoding/variance reduction to design a compressed DSGD algorithm.
\cite{Ma2021Springer, Xiong2022TAC, Liao2022TAC} developed quantized/compressed DGT algorithms by integrating the gradient tracking method with a uniform quantizer and compressors with bounded relative compression error compressors, respectively.
\cite{Yuan2012distributed} proposed a distributed averaging method using random quantization.

Previous studies have focused on distributed convex optimization.
However, many applications such as distributed learning \cite{omidshafiei2017PMLR}, distributed clustering \cite{forero2011JSTSP}, the cost functions are usually nonconvex, see, e.g., \cite{Hong2017PMLR,Yi2021TAC}.
Consequently, several studies have proposed compressed distributed nonconvex optimization algorithms.
For example, \cite{Taheri2020ICML,Reisizadeh2019NIPS} proposed compressed DSGD algorithms that utilize exact and random quantization, respectively.
\cite{xu2022Arxiv} proposed two quantized distributed algorithms by integrating a uniform quantizer with the DGT and DPI algorithms\TY{, respectively}.
\blue{\cite{Liao2023Arxiv} proposed a compressed DGT algorithm that utilizes robust compressors.}
\red{Furthermore, \cite{yi2022TAC} developed compressed distributed primal--dual algorithms, which used several general compressors.
In this paper, we investigate another distributed optimization algorithm equipped with \blue{these} general compressors, i.e., DGT algorithm.
As explained by \cite{Koloskova2021NIPS}, EXTRA \cite{Shi2015SIAM} and distributed primal--dual \cite{Alghunaim2020Aut} algorithms typically demand a noiseless setting.
In contrast, the DGT algorithm \blue{is able} to tolerate stochastic noise \cite{Di2016TSIPN,Nedic2017SIAM}, rendering it suitable for solving nonconvex optimization problems, especially in machine learning, see, e.g., \cite{Lin2021ICML,Yuan2021ICCV}.}
This motivates us to consider the gradient tracking methods with communication compression for the distributed nonconvex optimization problem.

The main contributions of this paper are:

\begin{itemize}
\item \blue{For the compressors with bounded relative compression error, which include the commonly considered unbiased and biased but contractive compressors,
we design a compressed DGT algorithm (Algorithm~\ref{thm-1}).
For smooth local cost functions, we \TY{design} an appropriate Lyapunov function \TY{in Theorem~\ref{thm-1}} to \TY{show} that the proposed algorithm sublinearly converges to a stationary
point.}
Moreover, if the global cost function satisfies the Polyak--{\L}ojasiewicz (P--{\L}) condition, which is weaker than the standard strong convexity condition and the global minimizer is not
necessarily unique, we establish \TY{in Theorem~\ref{thm-2}} that the proposed algorithm linearly converges to a global optimal point.

\item We propose an error feedback based compressed gradient tracking algorithm (Algorithm~\ref{alg:2}) to improve the algorithm's efficiency for biased compression methods.
Moreover, we redesign a Lyapunov function \TY{in Theorem~\ref{thm-3}} to accommodate the introduced error feedback variables, and utilize it to establish convergence results \TY{without
(Theorem~\ref{thm-3}) and with (Theorem~\ref{thm-4})} the P--{\L} condition.

\item \blue{For the compressors with bounded absolute compression error, which includes the commonly considered unbiased compressors with bounded variance, we develop a compressed DGT
algorithm (Algorithm~\ref{thm-3}), \TY{and redesign a Lyapunov function \TY{in Theorem~\ref{thm-5}}.}}
For compressors with globally bounded absolute compression error, we present the convergence results \TY{without (Theorem~\ref{thm-5}) and with (Theorem~\ref{thm-6})} the P--{\L}
condition, which are similar to Theorems~\ref{thm-1} and ~\ref{thm-2}.

\item \TY{F}or compressors with locally bounded absolute compression error, \TY{we redesign a Lyapunov function \TY{in Theorem~\ref{thm-7}}, to
demonstrate that the proposed Algorithm~\ref{thm-3}} linearly converge to a global optimal point under the P--{\L} condition.
\end{itemize}
\blue{Note that for perfect communication in DGT nonconvex optimization algorithms under the P--{\L} condition, \cite{Xin2021TSP, tang2020distributed} constructed systems of linear inequalities to analyze the convergence of the algorithms.
Nonetheless, \TY{the} prior knowledge is \TY{required} to determine \TY{their proposed algorithms'} parameters.
To avoid the need for the P--{\L} constant in determining algorithm parameters, this paper employs the Lyapunov method to analyze the convergence of the proposed algorithms.
To the best of our knowledge, this paper is the first to avoid using the P--{\L} constant in the context of the DGT nonconvex optimization algorithm under the P--{\L} condition.}
\TY{This is a significant property since determining the P--{\L} constant can be a challenging task.}
\red{Moreover, for the \blue{perfect} communication scenario, Laypunov analysis has been developed for DGT in both the convex \cite{Notarnicola2023TAC} and nonconvex cases \cite{Carnevale2022TCSL}.
The Lyapunov analysis method is simpler than other methods, which typically require constructing systems of linear inequalities, employing Lyapunov-like arguments, or utilizing control tools, see, e.g., \cite{Xin2021TSP, tang2020distributed,Varagnolo2015TAC,Qu2017TCNS,Xu2017TAC}.}
\lx{In \cite{zhao2022NIPS}, the authors also utilized the Lyapunov method to analyze the convergence of the proposed BEER algorithm.
However, the compressors with bounded relative compression error that we considered are more general than the biased but contractive compressors used in \cite{zhao2022NIPS}.
Furthermore, we also consider compressors with globally and locally bounded absolute compression error.}
\TY{To the best of our knowledge, this paper is the first to utilize the Lyapunov method for analyzing the compressed DGT algorithm in the context of nonconvex optimization problems.}

This work differs from \cite{yi2022TAC} in \xul{two} key aspects:
(i) This paper is built upon the DGT algorithm to design compressed distributed nonconvex optimization algorithms, whereas \cite{yi2022TAC} is based on the distributed primal--dual algorithm. The primary distinction between these compressed algorithms lies in the fact that the gradient tracking algorithm is better suited for the nonconvex setting. In the numerical simulation section, we showcase that the proposed algorithms achieve faster convergence when compared to the compressed algorithms \TY{proposed} in \cite{yi2022TAC}.
(ii) In this paper, for each compressed algorithm, we developed a novel Lyapunov function to analyze the convergence of the proposed algorithm. This Lyapunov function differs from the one \blue{proposed} in \cite{yi2022TAC}.

The remainder of the paper is organized as follows.
Section~\ref{sec-2} presents the problem formulation.
\TY{In Sections~\ref{sec-3} and~\ref{sec-4}, we introduce three compressed gradient tracking algorithms and conduct analyses for compressors involving bounded relative compression errors, as well as globally and locally bounded absolute compression errors, respectively.}
Section~\ref{sec-5} presents numerical simulation examples.
Finally, concluding remarks are offered in Section~\ref{sec-6}.

\textbf{Notation:} Let $\mathbf{1}_n$ (or $\mathbf{0}_n$) be the $n\times 1$ vector with all ones (or zeros), and $\mathbf{I}_n$ be the $n$-dimensional identity matrix.
$\mathrm{col}(Z_{1},\ldots,Z_{n})$ is the concatenated column vector of vectors $Z_{i}\in\mathbb{R}^{d}$.
$\|\cdot\|$ is the Euclidean vector norm or spectral matrix norm. For a column vector $X=(X_{1},\dots,X_{m})$,  $\|X\|_{\infty}=\max_{1\leq i\leq m}|X_{i}|$.
For a positive semi-definite matrix $\mathcal{M}$,
$\rho(\mathcal{M})$ is the spectral radius.
The minimum integer less than or equal to $c$ is denoted by $\lfloor c\rfloor$.
$\operatorname{sign}(c)$ and $|c|$ are the element--wise sign and absolute value, respectively.
\XL{Given any differentiable function $F$, $\nabla F$ is the gradient of $F$.}
$A\otimes B$ represents the Kronecker product of matrices $A$ and $B$.
$A\preceq B$ if all entries of matrix $A-B$ are not greater than zero, and $A\succ0$ if all entries of matrix $A$ that are greater than zero.
\xul{$\mathbb{Z}^{+}$ denotes the set of positive integers.}

\section{Problem Formulation}\label{sec-2}
Consider a group of $n$ agents distributed over a directed graph $\mathcal{G}=(\mathcal{V},\mathcal{E})$, where $\mathcal{V}=\{1,2,\ldots,n\}$ is the vertex set and $\mathcal{E}\subseteq\mathcal{V}\times\mathcal{V}$ is the set of directed edges. A directed path from agent $i_{1}$ to agent $i_{k}$ is a sequence of agents $\{i_{1},\dots,i_{k}\}$ such that $(i_{j},i_{j+1})\in\mathcal{E}$, \xul{$j=1,\dots,k-1$}. A directed graph is strongly connected if there exists a path between any pair of distinct agents.

Assume that each agent has a private differentiable local cost function $F_{i}:\mathbb{R}^d\rightarrow\mathbb{R}$, the optimal set $\mathbb{X}^{\star}=\emph{\emph{argmin}}_{X\in\mathbb{R}^{d}}F(X)$ is nonempty and $F^{\star}=\emph{\emph{min}}_{X\in\mathbb{R}^{d}}F(X)>-\infty$.
\XL{The objective is to find an optimizer $X^\star$ to minimize the average of all local cost functions $F(X) = \frac{1}{n} \sum_{i=1}^{n}F_{i}(X)$, that is,
\vspace{-3mm}
\begin{equation}
\min_{X \in \mathbb{R}^d}F(X)= \min_{X \in \mathbb{R}^d}\frac{1}{n} \sum_{i=1}^{n}F_{i}(X).
\end{equation}}

Throughout this paper, we make the following assumptions.

\begin{assumption}\label{assum-4}
The directed graph $\mathcal{G}$ is strongly connected and permits a nonnegative doubly stochastic weight matrix $W=[w_{ij}]\in\mathbb{R}^{n\times n}$, where $w_{ii}>0$, for all $i\in\mathcal{V}$, and $w_{ij}>0$ if and only if agent $i$ can receive information from agent $j$, otherwise
$w_{ij}=0$. Moreover, $W\bm{1}_{n}=\bm{1}_{n}$ and $\bm{1}_{n}^{T}W=\bm{1}_{n}^{T}$, ensuring that $W$ is doubly stochastic.
\end{assumption}


\begin{assumption}\label{assum-6}
Each local cost function $F_{i}(X)$ is smooth with constant $L_{f}>0$, i.e.,
\begin{equation}\label{assum2}
\|\nabla F_{i}(X)-\nabla F_{i}(Y)\|\leq L_{f}\|X-Y\|,~\forall X,Y\in\lx{\mathbb{R}^{d}.}
\end{equation}
\end{assumption}

\begin{assumption}\label{assum-7}
The global cost function $F(X)$ satisfies the Polyak--${\L}$ojasiewicz (P--{\L}) condition with constant $\nu>0$, i.e.,
  \begin{equation}\label{PL}
  \frac{1}{2}\|\nabla F(X)\|^{2}\geq\nu(F(X)-F^{\star}),~\forall X\in\mathbb{R}^{d}.
  \end{equation}
\end{assumption}

Assumptions \ref{assum-4} and \ref{assum-6} are common in the literature, e.g., \cite{nedic2018distributed,Yang2019ARC}.
\xul{Assumption 3 does not imply convexity of the global cost function, but it ensures that all stationary points are global optima.}

%

Motivated by scenarios where the communication channel often has limited bandwidth,
\lx{we propose three compressed distributed nonconvex (without and with the P--{\L} condition) optimization algorithms utilizing compressors with bounded relative compression error (Assumption~\ref{assum-1}), globally bounded absolute compression error (Assumption~\ref{assum-2}), and locally bounded absolute compression error (Assumption~\ref{assum-3}), respectively, in the subsequent \xul{sections}.}

\section{Compressed Distributed Nonconvex Algorithms: Bounded Relative Compression Error}\label{sec-3}
In this section, we introduce a compression operator with bounded relative compression error.
In Section~\ref{Relative-A}, we propose a compressed DGT algorithm, and present the convergence results. Additionally, in Section~\ref{Relative-B}, we extend the compressed algorithm to an error feedback version for biased compressors, and present the convergence result.
\begin{assumption}\label{assum-1}
\cite{yi2022TAC,Liao2022TAC} The compression operator $\mathcal{C}: \mathbb{R}^{d}\rightarrow\mathbb{R}^{d}$, adheres to the condition:
\begin{equation}\label{Assum1b}
\bm{\mathrm{E}}_{\mathcal{C}}[\|\frac{\mathcal{C}(X)}{r}-X\|^{2}]\leq(1-\psi)\|X\|^{2},~\forall X\in\mathbb{R}^{d},
\end{equation}
for some constants $r>0$ and $\psi\in(0,1]$. Here $\mathrm{E}_{\mathcal{C}}[\cdot]$ represents the expectation over the internal randomness of the compression operator $\mathcal{C}$.
\end{assumption}
This condition implies:
\begin{equation}\label{Assum1a}
\bm{\mathrm{E}}_{\mathcal{C}}[\|\mathcal{C}(X)-X\|^{2}]\leq C\|X\|^{2},~\forall X\in\mathbb{R}^{d},
\end{equation}
where $C=2r^{2}(1-\psi)+2(1-r)^{2}$.


Assumption~\ref{assum-1} encompasses various compression operators commonly used in the literature, such as norm-sign compression operators, random quantization, and sparsification \cite{yi2022TAC,Liao2022TAC}. It represents a broader class of compressors utilized in distributed optimization algorithms.



\subsection{Compressed Gradient Tracking Algorithm: Bounded Relative Compression Error}\label{Relative-A}

\TY{In this subsection, we propose the compressed DGT algorithm (Algorithm~\ref{alg:1}), which is the same as the C-GT proposed in \cite{Liao2022TAC}.
In \cite{Liao2022TAC}, the authors constructed systems of linear inequalities to demonstrate the linear convergence of the proposed algorithm under the strongly convex case.
In this subsection, \lx{we demonstrate that the proposed Algorithm~\ref{alg:1} exhibits sublinear convergence in the general nonconvex case (Theorem~\ref{thm-1}) and linear convergence when the P--{\L} condition is satisfied (Theorem~\ref{thm-2}).}
Furthermore, we design a Lyapunov function for analyzing the convergence of the proposed algorithm.
\lx{Benefiting from this Lyapunov function, in Theorem~\ref{thm-2}, we are able to design the algorithm parameters without prior knowledge of the P--{\L} constant}. This features a unique aspect that separates our work from the existing DGT nonconvex optimization results (\cite{Xin2021TSP, tang2020distributed}) under the P--{\L} condition.}

\begin{algorithm}[!h]
\caption{}
  \label{alg:1}
  \begin{algorithmic}
  \STATE  $\bm{\mathrm{For~each~agent}}$ $i\in\mathcal{V}.$
	\STATE $\bm{\mathrm{Initialization\hspace{-1mm}:}}$\\
     $X_{i}(0)\in\mathbb{R}^{d}$, $Y_{i}(0)=\nabla F_{i}(X_{i}(0))$, $A_{i}(0)=B_{i}(0)=C_{i}(0)=D_{i}(0)=\bm{0}_{d}$, $Q_{i}^{X}(0)=\mathcal{C}(X_{i}(0))$, and $Q_{i}^{Y}(0)=\mathcal{C}(Y_{i}(0))$.
\STATE $\bm{\mathrm{Communication\hspace{-1mm}:}}$\\
Transmit $Q_{i}^{X}(k)$ and $Q_{i}^{Y}(k)$ to its \xul{out-neighbors} and receive $Q_{j}^{X}(k)$ and $Q_{j}^{Y}(k)$ from its \xul{in-neighbors}.
	\STATE $\bm{\mathrm{Update~Rule\hspace{-1mm}:}}$
\begin{subequations}
\begin{align}
&\hspace{-2mm}A_{i}(k+1)=A_{i}(k)+\varphi_{X}Q_{i}^{X}(k),\label{Alg1b}\\
  &\hspace{-2mm}B_{i}(k+1)=B_{i}(k)+\varphi_{X}(Q_{i}^{X}(k)-\sum_{j=1}^{n}W_{ij}Q_{j}^{X}(k)),\label{Alg1c}\\
  &\hspace{-2mm}C_{i}(k+1)=C_{i}(k)+\varphi_{Y}Q_{i}^{Y}(k),\label{Alg1e}\\
  &\hspace{-2mm}D_{i}(k+1)=D_{i}(k)+\varphi_{Y}(Q_{i}^{Y}(k)-\sum_{j=1}^{n}W_{ij}Q_{j}^{Y}(k)),\label{Alg1f}\\
  &\hspace{-2mm}X_{i}(k+1)=X_{i}(k)-\gamma[B_{i}(k)+Q_{i}^{X}(k)-\sum_{j=1}^{n}W_{ij}Q_{j}^{X}(k)]-\eta Y_{i}(k),\label{Alg1g}\\
  &\hspace{-2mm}Y_{i}(k+1)=Y_{i}(k)-\gamma[D_{i}(k)+Q_{i}^{Y}(k)-\sum_{j=1}^{n}W_{ij}Q_{j}^{Y}(k)]\notag\\
  &\qquad\qquad\quad+\nabla F_{i}(X_{i}(k+1))-\nabla F_{i}(X_{i}(k)),\label{Alg1h}\\
  &\hspace{-2mm}Q_{i}^{X}(k+1)=\mathcal{C}(X_{i}(k+1)-A_{i}(k+1)),\label{Alg1a}\\
  &\hspace{-2mm}Q_{i}^{Y}(k+1)=\mathcal{C}(Y_{i}(k+1)-C_{i}(k+1)),\label{Alg1d}
  \end{align}
\end{subequations}
where $\gamma$, $\eta$, $\varphi_{X}$, and $\varphi_{Y}$ are positive parameters.
\end{algorithmic}
\end{algorithm}

\begin{remark}\label{rem-new}
\XL{The main difference between the proposed Algorithm~1 and BEER proposed in \cite{zhao2022NIPS} lies in the introduction of $\varphi_{X}$ and $\varphi_{Y}$ to correct the error caused by the $r$-scaling of the compression operator.
It is easy to verify that Algorithm 1 reduces to BEER when $\varphi_{X}=\varphi_{Y}=1$.
Due to the introduction of $\varphi_{X}$ and $\varphi_{Y}$, Algorithm~1 is suitable for a broader range of compression operators, including (\ref{Assum1b}).
Note that the compressors in (4) are equivalent to the compressors used in \cite{zhao2022NIPS} when $r=1$.
However, the introduction of $\varphi_{X}$ and $\varphi_{Y}$, along with $r\neq1$, makes it challenging to straightforwardly extend the proof techniques used in BEER.
More specifically, when analyzing the convergence of compression errors, we must construct inequalities that satisfy the properties of the compressors we consider. This requires frequent use of inequality shrinking techniques, and if the increased conservativeness from these shrinkings is not carefully managed, it may prevent us from proving algorithm convergence. These factors present significant challenges in our analysis.}
\end{remark}
We denote $\bm{X}=\mathrm{col}(X_{1},\ldots,X_{n})$, $\bm{Y}=\mathrm{col}(Y_{1},\ldots,Y_{n})$, $\bm{A}=\mathrm{col}(A_{1},\ldots,A_{n})$, $\bm{C}=\mathrm{col}(C_{1},\ldots,C_{n})$, $\bar{X}(k)=\frac{1}{n}(\bm{1}_{n}^{T}\otimes\bm{\mathrm{I}}_{d})\bm{X}(k)$, $\bar{\bm{X}}(k)=\bm{1}_{n}\otimes\bar{X}(k)$, $\bm{H}=\frac{1}{n}(\bm{1}_{n}\bm{1}^{T}_{n}\otimes\bm{\mathrm{I}}_{d})$, $\bar{\bm{Y}}(k)=\bm{H}\bm{Y}(k)$.

\blue{To analyze the convergence of Algorithm~\ref{alg:1}, we consider the following Lyapunov candidate function.}

\begin{equation}\label{Lyapunov}
U(k)=V(k)+n(F(\bar{X}(k))-F^{\star}),
\end{equation}
where
$$
V(k)=\|\bm{X}(k)-\bar{\bm{X}}(k)\|^{2}+\phi\|\bm{Y}(k)-\bar{\bm{Y}}(k)\|^{2}+\|\bm{X}(k)-\bm{A}(k)\|^{2}+\|\bm{Y}(k)-\bm{C}(k)\|^{2},
$$
and $\phi=\frac{(1-\sigma)^{2}}{320L_{f}^{2}}$ \lx{with $\sigma\in(0,1)$ being the spectral norm of $W-\frac{1}{n}\bm{1}_{n}\bm{1}^{T}_{n}$.}

Note that the designed \red{Lyapunov candidate} function (\ref{Lyapunov}) incorporates several error terms: \lx{the} consensus error term $\|\bm{X}(k)-\bar{\bm{X}}(k)\|^{2}$, gradient tracking error term $\|\bm{Y}(k)-\bar{\bm{Y}}(k)\|^{2}$, compression error terms $\|\bm{X}(k)-\bm{A}(k)\|^{2}$ and $\|\bm{Y}(k)-\bm{C}(k)\|^{2}$, and the optimal error term $n(F(\bar{X}(k))-F^{\star})$.
\TY{The weight parameter $\phi$ is instrumental in fine-tuning the values of the respective terms within the designed Lyapunov function, thereby ensuring the convergence of the proposed algorithms.}

\blue{We are now ready to present the convergence results of Algorithm~\ref{alg:1}.}


%
\begin{theorem}\label{thm-1}
Suppose that Assumptions~\ref{assum-4}--\ref{assum-6}, and \ref{assum-1} hold.
Let each agent $i\in\mathcal{V}$ run Algorithm~\ref{alg:1} with algorithm parameters $\varphi_{X}$, $\varphi_{Y}\in(0,\frac{1}{r})$, $\eta$ and $\gamma$ such that
\begin{align}\label{eta}
\eta\in&(0,\min\{\frac{(1-\sigma)^{2}\gamma}{40L_{f}},~\frac{0.4(1-\sigma)\gamma}{L_{f}^{2}},~\frac{(1-\sigma)^{2}}{80L_{f}}\sqrt{\frac{\gamma}{1+c_{1}^{-1}}},\notag\\
&\frac{9}{40(4(1+c_{1}^{-1})+5(1+c_{2}^{-1}))},~\frac{1}{2L_{f}},~\gamma\}),
\end{align}
\begin{align}\label{gamma}
\gamma&\in(0,\Pi:=\min\{\frac{1-\sigma}{160(1+c_{1}^{-1})}, ~\frac{1-\sigma}{40000(1+c_{2}^{-1})L_{f}^{2}},\notag\\
&\frac{c_{1}(1-\sigma)}{40C},~\frac{c_{1}}{8\sqrt{C}},~\frac{c_{1}}{10L_{f}\sqrt{C(1+c_{2}^{-1})}},\frac{c_{2}L_{f}^{2}}{C},~\frac{c_{2}}{10\sqrt{C}}\}),
\end{align}
where $c_{1}=\frac{\varphi_{X}\psi r}{2}$,
$c_{2}=\frac{\varphi_{Y}\psi r}{2}$, and $C=2r^{2}(1-\psi)+2(1-r)^{2}$.

Then, we have
\begin{equation}\label{thm-1a}
\sum_{t=0}^{k}\bm{\mathrm{E}}_{\mathcal{C}}[\|\bm{X}(t)-\bar{\bm{X}}(t)\|^{2}+n\|\nabla F(\bar{X}(k))\|^{2}]\leq\frac{U(0)}{\theta_{1}},
\end{equation}
and
\begin{equation}\label{thm-1b}
\bm{\mathrm{E}}_{\mathcal{C}}[n(F(\bar{X}(k))-F^{\star})]<U(0),
\end{equation}
where
\begin{align*}
\theta_{1}&=\min\{\theta_{2},\theta_{3}\},~\theta_{2}=\frac{\eta}{4}-(\phi\varepsilon_{1}+\varepsilon_{2}+\varepsilon_{3}),\\
\theta_{3}&=\min\{0.07(1-\sigma)\gamma,0.44c_{1}(2c_{1}+1),0.77c_{2}(2c_{2}+1)\},\\
\varepsilon_{1}&=\frac{8L_{f}^{2}}{(1-\sigma)\gamma}\eta^{2},~\varepsilon_{2}=4(1+c_{1}^{-1})\eta^{2},~\varepsilon_{3}=5(1+c_{2}^{-1})\eta^{2}.
\end{align*}
\end{theorem}
\begin{IEEEproof}
The proof is given in Appendix~B.
\end{IEEEproof}

\begin{remark}\label{rem-4}
Theorem~\ref{thm-1} shows that Algorithm~\ref{alg:1} achieves a convergence rate of $\mathcal{O}(1/T)$. Specifically, (\ref{thm-1a}) reveals that the term $\min_{k\leq T}\{\bm{\mathrm{E}}_{\mathcal{C}}[n\|\nabla F(\bar{X}(k))\|^2+\|\bm{X}(k)-\bar{\bm{X}}(k)\|^2]\}$ decays at a rate of $\mathcal{O}(1/T)$. Additionally, (\ref{thm-1b}) indicates that the term $\bm{\mathrm{E}}_{\mathcal{C}}[n(F(\bar{X}(k))-F^{\star})]$ is bounded.
%
%
\end{remark}

\blue{Moreover, with Assumption~\ref{assum-7}, the following result shows that Algorithm~\ref{alg:1} can find global optima and the convergence rate is linear.}

\begin{theorem}\label{thm-2}
Suppose that Assumptions~\ref{assum-4}--\ref{assum-1} hold. Let each agent $i\in\mathcal{V}$ run Algorithm~\ref{alg:1} with the parameters $\eta$, $\gamma$, $\varphi_{X}$, and $\varphi_{Y}$ being given in Theorem~\ref{thm-1}.
Then, we have
\begin{equation}\label{thm-2a}
\bm{\mathrm{E}}_{\mathcal{C}}[\|\bm{X}(k)-\bar{\bm{X}}(k)\|^{2}+n(F(\bar{X}(k))-F^{\star})]\leq(1-\theta_{4})^{k}U(0),
\end{equation}
where $\theta_{4}=\min\{\theta_{3},2\nu\theta_{2}\}$.
\end{theorem}
\begin{IEEEproof}
The proof is given in Appendix~C.
\end{IEEEproof}

\begin{remark}\label{rem-2}
\TY{
Note that the P--{\L} constant is not utilized in Algorithm~\ref{alg:1}.
This is a significant property since determining the P--{\L} constant can be a challenging task.
\xul{It is worth noting that most existing DGT nonconvex optimization algorithms require the use of the P--{\L} constant}, see, e.g., \cite{Xin2021TSP,tang2020distributed,Liao2023Arxiv}.
This property of not requiring the P--{\L} constant arises from the Lyapunov method, which differs from methods based on constructing systems of linear inequalities, as used in \cite{Xin2021TSP, tang2020distributed, Liao2023Arxiv}.
}
\end{remark}

\subsection{Error Feedback Based Compressed Gradient Tracking Algorithm: Bounded Relative Compression Error}\label{Relative-B}
\TY{In this subsection, we extend Algorithm~\ref{alg:1} to an error feedback version for biased compressors, as shown in Algorithm~\ref{alg:2}, which is the same as the EF-C-GT proposed in \cite{Liao2022TAC}.
Similar to Section~\ref{Relative-A}, we investigate the nonconvex optimization problem without (Theorem~\ref{thm-3}) and with (Theorem~\ref{thm-4}) the P--{\L} condition. In this subsection, we reconstruct a Lyapunov function to analyze the convergence of Algorithm~\ref{alg:2}.}
\begin{algorithm}[htp]
\caption{}
  \label{alg:2}
  \begin{algorithmic}
  \STATE  $\bm{\mathrm{For~each~agent}}$ $i\in\mathcal{V}.$
	\STATE $\bm{\mathrm{Initialization:\hspace{-1mm}}}$\\
     $X_{i}(0)\in\mathbb{R}^{d}$, $Y_{i}(0)=\nabla F_{i}(X_{i}(0))$, $A_{i}(0)=B_{i}(0)=C_{i}(0)=D_{i}(0)=\bm{0}_{d}$, $Q_{i}^{X}(0)=\hat{Q}_{i}^{X}(0)=\mathcal{C}(X_{i}(0))$, and $Q_{i}^{Y}(0)=\hat{Q}_{i}^{Y}(0)=\mathcal{C}(Y_{i}(0))$.
\STATE $\bm{\mathrm{Communication:\hspace{-1mm}}}$\\
Transmit $Q_{i}^{X}(k)$, $\hat{Q}_{i}^{X}(k)$, $Q_{i}^{Y}(k)$, and $\hat{Q}_{i}^{Y}(k)$ to its \xul{out-neighbors} and receive $Q_{j}^{X}(k)$, $\hat{Q}_{j}^{X}(k)$, $Q_{j}^{Y}(k)$, and $\hat{Q}_{j}^{Y}(k)$ from its \xul{in-neighbors}.
	\STATE $\bm{\mathrm{Update~Rule:\hspace{-1mm}}}$
\hspace{-1cm}
\begin{subequations}
\begin{align}
&A_{i}(k+1)=A_{i}(k)+\varphi_{X}Q_{i}^{X}(k),\label{Alg2b}\\
&B_{i}(k+1)=B_{i}(k)+\varphi_{X}(Q_{i}^{X}(k)-\sum_{j=1}^{n}W_{ij}Q_{j}^{X}(k)),\label{Alg2c}\\
&C_{i}(k+1)=C_{i}(k)+\varphi_{Y}Q_{i}^{Y}(k),\label{Alg2e}\\
&D_{i}(k+1)=D_{i}(k)+\varphi_{Y}(Q_{i}^{Y}(k)-\sum_{j=1}^{n}W_{ij}Q_{j}^{Y}(k)),\label{Alg2f}\\
&X_{i}(k+1)=X_{i}(k)-\gamma[B_{i}(k)+\hat{Q}_{i}^{X}(k)-\sum_{j=1}^{n}W_{ij}\hat{Q}_{j}^{X}(k)]-\eta Y_{i}(k),\label{Alg2g}\\
&Y_{i}(k+1)=Y_{i}(k)-\gamma[D_{i}(k)+\hat{Q}_{i}^{Y}(k)-\sum_{j=1}^{n}W_{ij}\hat{Q}_{j}^{Y}(k)]\notag\\
&\qquad\qquad\quad+\nabla F_{i}(X_{i}(k+1))-\nabla F_{i}(X_{i}(k)),\label{Alg2h}\\
&Q_{i}^{X}(k+1)=\mathcal{C}(X_{i}(k+1)-A_{i}(k+1)),\label{Alg2a}\\
&Q_{i}^{Y}(k+1)=\mathcal{C}(Y_{i}(k+1)-C_{i}(k+1)),\label{Alg2d}\\
&E_{i}^{X}(k+1)=\varsigma E_{i}^{X}(k)+X_{i}(k)-A_{i}(k)-\hat{Q}_{i}^{X}(k),\label{Alg2i}\\
&E_{i}^{Y}(k+1)=\varsigma E_{i}^{Y}(k)+Y_{i}(k)-C_{i}(k)-\hat{Q}_{i}^{Y}(k),\label{Alg2j}\\
&\hat{Q}_{i}^{X}(k+1)=\mathcal{C}(\varsigma E_{i}^{X}(k+1)+X_{i}(k+1)-A_{i}(k+1)),\label{Alg2k}\\
&\hat{Q}_{i}^{Y}(k+1)=\mathcal{C}(\varsigma E_{i}^{Y}(k+1)+Y_{i}(k+1)-C_{i}(k+1)),\label{Alg2l}
  \end{align}
\end{subequations}
where $\gamma$, $\eta$, $\varphi_{X}$, $\varphi_{Y}$, and $\varsigma$ are positive parameters.
\end{algorithmic}
\end{algorithm}

Before demonstrating the convergence of Algorithm~\ref{alg:2},
we denote
$\bm{E}^{X}=\mathrm{col}(E^{X}_{1},\ldots,E^{X}_{n})$,
$\bm{E}^{Y}=\mathrm{col}(E^{Y}_{1},\ldots,E^{Y}_{n})$,
\blue{To analyze the convergence of Algorithm~\ref{alg:1}, we consider the following Lyapunov candidate function.}
\begin{equation}\label{Lyapunov2}
\hat{U}(k)=\hat{V}(k)+n(F(\bar{X}(k))-F^{\star}),
\end{equation}
where
$$
\hat{V}(k)=V(k)+\hat{\phi}(\|\bm{E}^{X}(k)\|^{2}+\|\bm{E}^{Y}(k)\|^{2}),~\hat{\phi}=\frac{0.1}{C}\min\{c_{1}(2c_{1}+1),c_{2}(2c_{2}+1)\}.
$$

Note that Algorithm~\ref{alg:2} introduces two error feedback variables, $\bm{E}^{X}$ and $\bm{E}^{Y}$, to rectify the bias caused by the biased compressors.
Hence, the designed \red{Lyapunov candidate} function (\ref{Lyapunov2}) incorporates two additional feedback error terms, $\|\bm{E}^{X}(k)\|^{2}$ and \blue{$\|\bm{E}^{Y}(k)\|^{2}$}.
Moreover, the weight parameter $\hat{\phi}$ plays a crucial role in ensuring the convergence of this function.

Next, we investigate the convergence of Algorithm~\ref{alg:2}.
Similar to Theorem~\ref{thm-1}, we first establish  the following sublinear convergence result for Algorithm \ref{alg:2} without the P--{\L} condition.

\begin{theorem}\label{thm-3}
Suppose that Assumptions~\ref{assum-4}--\ref{assum-6}, and \ref{assum-1} hold.
Let each agent $i\in\mathcal{V}$ run Algorithm~\ref{alg:2} with parameters $\eta$, $\varphi_{X}$, and $\varphi_{Y}$ being chosen in Theorem~\ref{thm-1}, and $\gamma$, $\varsigma$, satisfying
\begin{align}\label{eta2}
\gamma\in&(0,~\min\{\frac{c_{1}(1-\sigma)}{160C},~\frac{c_{1}}{16\sqrt{C}},~\frac{c_{1}}{20L_{f}\sqrt{C(1+c_{2}^{-1})}},~\frac{c_{2}L_{f}^{2}}{4C},~\frac{c_{2}}{20\sqrt{C}},\notag\\
&\frac{1}{4(1+c_{1}^{-1})+5(1+c_{2}^{-1})L_{f}^{2}},~\frac{(1-\sigma)\hat{\phi}}{4(16(1+4\phi L_{f}^{2})+8(1-\sigma))},\notag\\
&\frac{1}{5(1+c_{2}^{-1})},~\frac{(1-\sigma)\hat{\phi}}{32(2\phi+(1-\sigma))},~\Pi\}),
\end{align}
\begin{equation}\label{vartheta}
\varsigma\in(0,~\min\{\frac{1}{2\sqrt{C}},~\frac{1}{\sqrt{2C+1}}\}).
\end{equation}
%

Then, we have
\begin{equation}\label{thm-3a}
\sum_{t=0}^{k}\bm{\mathrm{E}}_{\mathcal{C}}[\|\bm{X}(t)-\bar{\bm{X}}(t)\|^{2}+n\|\nabla F(\bar{X}(k))\|^{2}]\leq\frac{\hat{U}(0)}{\hat{\theta}_{1}},
\end{equation}
and
\begin{equation}\label{thm-3b}
\bm{\mathrm{E}}_{\mathcal{C}}[n(F(\bar{X}(k))-F^{\star})]\leq\bm{\mathrm{E}}_{\mathcal{C}}[\hat{U}(k)]<\hat{U}(0),
\end{equation}
where
$$
\hat{\theta}_{1}=\min\{\theta_{2},\hat{\theta}_{2}\},~\hat{\theta}_{2}=\min\{0.07(1-\sigma)\gamma,~0.24c_{1}(2c_{1}+1),~0.57c_{2}(2c_{2}+1),~0.25\}.
$$
\end{theorem}
\begin{IEEEproof}
The proof is given in Appendix~D.
\end{IEEEproof}


Similar to Theorem~\ref{thm-2}, we then have the following linear convergence result for Algorithm~\ref{alg:2} with Assumption~\ref{assum-7}.

\begin{theorem}\label{thm-4}
Suppose that Assumptions~\ref{assum-4}--\ref{assum-1} hold. Let each agent $i\in\mathcal{V}$ run Algorithm~\ref{alg:2} with
parameters $\eta$,
$\gamma$, $\varsigma$, $\varphi_{X}$, and $\varphi_{Y}$ being given in Theorem~\ref{thm-3}.
Then, we have
\begin{equation}\label{thm-4a}
\bm{\mathrm{E}}_{\mathcal{C}}[\|\bm{X}(k)-\bar{\bm{X}}(k)\|^{2}+n(F(\bar{X}(k))-F^{\star})]\leq(1-\hat{\theta}_{3})^{k}U(0),
\end{equation}
where $\hat{\theta}_{3}=\min\{\hat{\theta}_{2},2\nu\theta_{2}\}$.
\end{theorem}
\begin{IEEEproof}
The proof is given in Appendix~E.
\end{IEEEproof}

\section{Compressed Distributed Nonconvex Algorithm: Bounded Absolute Compression Error}\label{sec-4}
\TY{In this section, we propose a compressed DGT algorithm (Algorithm~\ref{alg:3}) that is designed for compressors with bounded absolute compression error, which is similar to the RCPP algorithm proposed in \cite{Liao2023Arxiv}.
In \cite{Liao2023Arxiv}, the authors constructed systems of linear inequalities to analyze the convergence of RCPP algorithm using more general compressors, \xul{which allow both locally and globally absolute compression errors,} under the P--{\L} condition for directed graphs.
\lx{However, \cite{Liao2023Arxiv} requires prior knowledge of the P--{\L} constant to design the algorithm parameters.}
\lx{In Section~\ref{Absolute-A}, we employ the Lyapunov method to demonstrate that the proposed Algorithm~\ref{alg:3} exhibits sublinear convergence in the general nonconvex case (Theorem~\ref{thm-5}) and linear convergence when the P--{\L} condition is satisfied (Theorem~\ref{thm-6}).
It is important to note that the P--{\L} constant is not used in designing the algorithm's parameters.}
In Section~\ref{Absolute-B}, we focus on compressors with locally bounded compression errors and establish a linear convergence result for Algorithm~\ref{alg:3} with P--{\L} condition \xul{(Theorem~\ref{thm-7})}.}

\begin{algorithm}[h!]
\caption{}
  \label{alg:3}
  \begin{algorithmic}
  \STATE  $\bm{\mathrm{For~each~agent}}$ $i\in\mathcal{V}.$
	\STATE $\bm{\mathrm{Initialization:\hspace{-1mm}}}$\\
     $X_{i}(0)\in\mathbb{R}^{d}$, $Y_{i}(0)=\nabla F_{i}(X_{i}(0))$, $\hat{X}_{i}(-1)=\hat{Y}_{i}(-1)=V_{i}(-1)=Z_{i}(-1)=\bm{0}_{d}$, $Q_{i}^{X}(0)=\mathcal{C}(X_{i}(0)/s(0))$, and $Q_{i}^{Y}(0)=\mathcal{C}(Y_{i}(0)/s(0))$, and $s(0)>0$.
\STATE $\bm{\mathrm{Communication:\hspace{-1mm}}}$\\
Transmit $Q_{i}^{X}(k)$ and $Q_{i}^{Y}(k)$ to its \xul{out-neighbors} and receive $Q_{j}^{X}(k)$ and $Q_{j}^{Y}(k)$ from its \xul{in-neighbors}.
	\STATE $\bm{\mathrm{Update~Rule:\hspace{-1mm}}}$
\begin{subequations}\label{Alg3}
\begin{align}
&\hat{X}_{i}(k)=\hat{X}_{i}(k-1)+s(k)Q_{i}^{X}(k),\label{Alg3a}\\
  &V_{i}(k)=V_{i}(k-1)+s(k)Q_{i}^{X}(k)-s(k)\sum_{j=1}^{n}W_{ij}Q_{j}^{X}(k),\label{Alg3b}\\
  &\hat{Y}_{i}(k)=\hat{Y}_{i}(k-1)+s(k)Q_{i}^{Y}(k),\label{Alg3c}\\
  &Z_{i}(k)=Z_{i}(k-1)+s(k)Q_{i}^{Y}(k)-s(k)\sum_{j=1}^{n}W_{ij}Q_{j}^{Y}(k),\label{Alg3d}\\
  &X_{i}(k+1)=X_{i}(k)-\gamma(\hat{X}_{i}(k)-V_{i}(k))-\eta Y_{i}(k),\label{Alg3e}\\
  &Y_{i}(k+1)=Y_{i}(k)-\gamma(\hat{Y}_{i}(k)-Z_{i}(k))+\nabla F_{i}(X_{i}(k+1))-\nabla F_{i}(X_{i}(k)),\label{Alg3f}\\
  &Q_{i}^{X}(k+1)=\mathcal{C}((X_{i}(k+1)-\hat{X}_{i}(k))/s(k+1)),\label{Alg3g}\\
  &Q_{i}^{Y}(k+1)=\mathcal{C}((Y_{i}(k+1)-\hat{Y}_{i}(k))/s(k+1)),\label{Alg3h}
  \end{align}
\end{subequations}
where $\gamma$, $\eta$, and $\mu$ are positive parameters, $s(k)=s(0)\mu^{k}>0$ is a decreasing scaling function, and $\mu\in(0,1)$.
\end{algorithmic}
\end{algorithm}

\subsection{Compressed Gradient Tracking Algorithm: Globally Bounded Absolute Compression Error}\label{Absolute-A}
In this subsection, we introduce a compression operator with globally bounded absolute compression error.
\begin{assumption}\label{assum-2}
\cite{yi2022TAC,Khirirat2020TSP}
The compression operator $\mathcal{C}: \mathbb{R}^{d}\rightarrow\mathbb{R}^{d}$, adheres to the condition:
\begin{equation}\label{Assum2}
\bm{\mathrm{E}}_{\mathcal{C}}[\|\mathcal{C}(X)-X\|_{p}^{2}]\leq C,~\forall X\in\mathbb{R}^{d},
\end{equation}
for $p\in\mathbb{Z}^{+}$ and constant $C\geq0$.
\end{assumption}
Assumption~\ref{assum-2} mainly includes deterministic quantization and unbiased random quantization. It is a commonly used compressor in the literature, as seen in \cite{yi2022TAC,Khirirat2020TSP}.
\blue{To analyze the convergence of Algorithm~\ref{alg:3} using compressors that have locally bounded absolute compression error, we consider the following Lyapunov candidate function:}
\begin{equation}\label{Lyapunov3}
\breve{U}(k)=\breve{V}(k)+n(F(\bar{X}(k))-F^{\star}),
\end{equation}
where
$$
\breve{V}(k)=\|\bm{X}(k)-\bar{\bm{X}}(k)\|^{2}+\phi\|\bm{Y}(k)-\bar{\bm{Y}}(k)\|^{2}.
$$

\xul{From (\ref{Alg3g}) and (\ref{Alg3h})}, it can be found that the compression errors are bounded by the scaling function $s(k)$.
Consequently, the reconstructed \red{Lyapunov candidate} function only comprises the consensus error term, gradient tracking error term, and the optimal error term.

Next, we investigate the convergence of Algorithm~\ref{alg:3}.

\begin{theorem}\label{thm-5}
Suppose that Assumptions~\ref{assum-4}--\ref{assum-6}, and \ref{assum-2} hold.
Let each agent $i\in\mathcal{V}$ run Algorithm~\ref{alg:3} with $s(0)>0$, and $\mu$ being an arbitrary constant in $(0,1)$,
parameters $\eta$ and $\gamma$ being chosen in Theorem~\ref{thm-1}.
Then, we have
\begin{equation}\label{thm-5a}
\sum_{t=0}^{k}\bm{\mathrm{E}}_{\mathcal{C}}[\|\bm{X}(t)-\bar{\bm{X}}(t)\|^{2}+n\|\nabla F(\bar{X}(k))\|^{2}]\leq\frac{\breve{U}(0)+\frac{\breve\theta_{2}}{1-\mu^{2}}}{\breve\theta_{1}},
\end{equation}
and
\begin{equation}\label{thm-5b}
\bm{\mathrm{E}}_{\mathcal{C}}[n(F(\bar{X}(k))-F^{\star})]\leq\bm{\mathrm{E}}_{\mathcal{C}}[\breve{U}(k)]<\breve{U}(0)+\frac{\breve\theta_{2}}{1-\mu^{2}},
\end{equation}
where
$$
\breve\theta_{1}=\min\{\breve\theta_{3},\breve\theta_{4}\},
~\xul{\breve\theta_{2}=\frac{16\gamma Cn\tilde{d}^2s^{2}(0)(1+2L_{f}^{2})}{1-\sigma}},~\breve\theta_{3}=0.59(1-\sigma)\gamma,~\breve\theta_{4}=\frac{\eta}{4}-\phi\varepsilon_{1},
$$
and \xul{$~\tilde{d}=1$ for $p\in[1,2]$, $\tilde{d}=d^{\frac{1}{2}-\frac{1}{p}}$ for $p>2$, with $p\in\mathbb{Z}^{+}$ representing the constant of the $p$-norm.}
\end{theorem}
\begin{IEEEproof}
The proof is given in Appendix~F.
\end{IEEEproof}

%
Similar to Theorem~\ref{thm-2}, we then have the following linear convergence result for Algorithm~\ref{alg:3} with Assumption~\ref{assum-7}.

\begin{theorem}\label{thm-6}
Suppose that Assumptions~\ref{assum-4}--\ref{assum-7}, and \ref{assum-2} hold.
Let each agent $i\in\mathcal{V}$ run Algorithm~\ref{alg:3} with $s(0)>0$, and $\mu$ being an arbitrary constant in $(0,1)$,
parameters $\eta$ and $\gamma$ being chosen in Theorem~\ref{thm-1}.
Then, we have
\begin{equation}\label{thm-6a}
\bm{\mathrm{E}}_{\mathcal{C}}[\|\bm{X}(k)-\bar{\bm{X}}(k)\|^{2}+n(F(\bar{X}(k))-F^{\star})]<(1-\breve\theta_{5})^{k}\breve\theta_{6},
\end{equation}
where
\begin{align*}
\breve\theta_{5}&=\min\{\breve\theta_{7},1-\mu^{2}\},\\
\breve\theta_{6}&=\breve{U}(0) + \breve\theta_{8}s^{2}(0)\begin{cases}
\frac{1}{(1-\breve\theta_{9})\mu^{2}}, & \text { if } 1-\breve\theta_{7}<\mu^{2}, \\
\frac{1}{(1-\breve\theta_{10})(1-\breve\theta_{7})}, & \text { if } 1-\breve\theta_{7}>\mu^{2}, \\
\frac{1}{(1-\breve\theta_{11})\varpi}, & \text { if } 1-\breve\theta_{7}=\mu^{2},
\end{cases}\\
\breve\theta_{7}&=\min\{\breve\theta_{1},2\nu\theta_{2}\},~\xul{\breve\theta_{8}=\frac{16n\tilde{d}^2\gamma C}{1-\sigma} (1+2L_{f}^{2})},\\
\breve\theta_{9}&=\frac{1-\breve\theta_{7}}{\mu^{2}},~\breve\theta_{10}=\frac{\mu^{2}}{1-\breve\theta_{7}},~\breve\theta_{11}=\frac{\mu^{2}}{\varpi},~\varpi\in(\mu^{2},1).
\end{align*}
\end{theorem}
\begin{IEEEproof}
The proof is given in Appendix~G.
\end{IEEEproof}


\subsection{Compressed Gradient Tracking Algorithm: Locally Bounded Absolute Compression Error}\label{Absolute-B}
In this subsection, we introduce a compression operator with locally bounded absolute compression error.
\begin{assumption}\label{assum-3}
\cite{yi2022TAC,Zhang2023TAC}
The compression operator $\mathcal{C}: \mathbb{R}^{d}\rightarrow\mathbb{R}^{d}$, adheres to the condition:
\begin{equation}\label{Assum3}
\|\mathcal{C}(X)-X\|_{p}\leq(1-\varphi),~\forall X\in\{X\in\mathbb{R}^{d}:\|X\|_{p}\leq1\},
\end{equation}
for $p\in\mathbb{Z}^{+}$ and constant $\varphi\in(0,1]$.
\end{assumption}

Assumption~\ref{assum-3} mainly includes standard quantization with dynamic and fixed quantization levels, which are commonly used compression techniques in the literature, see \cite{yi2022TAC,Zhang2023TAC}.

To analyze the convergence of Algorithm~\ref{alg:3} using compressors that have locally bounded absolute compression error, we consider the following Lyapunov candidate function:
\begin{equation}\label{Lyapunov4}
\tilde{U}(k)=\breve{V}(k)+\tilde{\phi}n(F(\bar{X}(k))-F^{\star}),
\end{equation}
where $\tilde{\phi}=\frac{0.4\gamma(1-\sigma)}{\eta L_{f}^{2}}$.

Next, we investigate the convergence of Algorithm~\ref{alg:3}.

\begin{theorem}\label{thm-7}
Suppose that Assumptions~\ref{assum-4}--\ref{assum-7} and \ref{assum-3} hold. Let each agent $i\in\mathcal{V}$ run Algorithm~\ref{alg:3} with the following given parameters:
\begin{equation}\label{eta4}
\eta\in(0,\min\{\frac{(1-\sigma)^{2}\gamma}{40L_{f}},~\frac{\varphi+\varphi^{2}-\varphi^{3}}{2\tilde{\xi}_{1}},~\frac{\varphi+\varphi^{2}-\varphi^{3}}{2\tilde{\xi}_{2}},1\}),
\end{equation}
\begin{equation}\label{gamma4}
\gamma\in(0,\min\{\sqrt{\frac{\varphi+\varphi^{2}-\varphi^{3}}{2\tilde{\xi}_{3}}},\sqrt{\frac{\varphi+\varphi^{2}-\varphi^{3}}{2\tilde{\xi}_{4}}},\frac{2L_{f}^{2}}{(1-\sigma)\nu}\}),
\end{equation}
\begin{equation}\label{s0}
s(0)\geq\max\{\sqrt{\frac{\tilde U(0)}{\tilde{\xi}_{5}}}, \max_{i\in\mathcal{V}}\|X_{i}(0)\|,\max_{i\in\mathcal{V}}\|Y_{i}(0)\|\},
\end{equation}
\begin{equation}\label{mu}
\mu\in[\max\{\sqrt{\tilde{\theta}_{1}},\sqrt{\tilde{\xi}_{6}},\sqrt{\tilde{\xi}_{7}}\},1),
\end{equation}
where
\begin{align*}
\tilde{\theta}_{1}&=1-\tilde{\theta}_{3}+\frac{\tilde{\theta}_{2}}{\tilde\xi_{5}}\gamma,~\tilde{\theta}_{2}=2(1+2L_{f}^{2})\frac{8n\tilde{d}^2(1-\varphi)^{2}}{1-\sigma},~\tilde{\theta}_{3}=\min\{\breve\theta_{3},0.5\gamma\},\\
\tilde{\theta}_{4}&=\min\{0.59(1-\sigma),~\frac{48\nu\phi}{1-\sigma}\},~\tilde{\xi}_{1}=4\hat{d}^{2}(5+4L_{f}^{2})(1+\varphi^{-1})\tilde{\xi}_{5},\\
\tilde{\xi}_{2}&=10L_{f}^{2}(3+2L_{f}^{2})(1+\varphi^{-1})\tilde{\xi}_{5},~\tilde{\xi}_{3}=\tilde{\xi}_{8}(1-\varphi)^{2}+32\hat{d}^{2}(1+\varphi^{-1})\tilde{\xi}_{5},\\
\tilde{\xi}_{4}&=\tilde{\xi}_{9}(1+L_{f}^{2})(1-\varphi)^{2}+40\hat{d}^{2}(1+L_{f}^{2})(1+\varphi^{-1})\tilde{\xi}_{5},\\
\tilde{\xi}_{5}&>\frac{\tilde{\theta}_{2}}{\tilde{\theta}_{4}},~\tilde{\xi}_{6}=(1-(\varphi+\varphi^{2}-\varphi^{3})+\eta\tilde{\xi}_{1}+\gamma^{2}\tilde{\xi}_{3}),\\
\tilde{\xi}_{7}&=(1-(\varphi+\varphi^{2}-\varphi^{3})+\eta\tilde{\xi}_{2}+\gamma^{2}\tilde{\xi}_{4}),\\
\tilde{\xi}_{8}&=16n\hat{d}^{2}\tilde{d}^{2}(1+\varphi^{-1}),~\tilde{\xi}_{9}=20n\hat{d}^{2}\tilde{d}^{2}(1+\varphi^{-1}),
\end{align*}
and \xul{$\hat{d}=d^{\frac{1}{2}-\frac{1}{p}}$ for $p\in[1,2]$, $\hat{d}=1$ for $p>2$, with $p\in\mathbb{Z}^{+}$ representing the constant of the $p$-norm.}
%
%

Then, we have
\begin{equation}\label{thm-7a}
\|\bm{X}(k)-\bar{\bm{X}}(k)\|^{2}+n(F(\bar{X}(k))-F^{\star})\leq\tilde{\xi}_{5}s^{2}(k).
\end{equation}
\end{theorem}
\begin{IEEEproof}
The proof is given in Appendix~H.
\end{IEEEproof}

\begin{remark}\label{rem-4}
\blue{The standard uniform quantizer, as used in \cite{Ma2021Springer, Xiong2022TAC} for the strongly convex case, and in \cite{xu2022Arxiv} for the nonconvex case under the P--{\L} condition, \xul{is a widely used method for reducing communication overhead in distributed optimization}. Note that Assumption~\ref{assum-3} serves as a more general compressor. In other words, Theorem~\ref{thm-7} demonstrates that the proposed algorithm achieves linear convergence with a broader range of compressors and only requires the global cost function to satisfy the P--{\L} condition.}
\end{remark}

\section{Numerical Examples}\label{sec-5}
In this section, we demonstrate the effectiveness of the proposed compressed distributed algorithms.
Consider \red{a strongly connected and weight-balanced directed network} consisting of $20$ agents and the communication graph is randomly generated.
Then, we consider the distributed nonconvex optimization problem studied in \cite{tang2020distributed}, i.e.,
  \begin{equation}\label{numerical}
F_{i}(X)=\frac{h_{i}}{({1+\exp\left(-\xi_i^{T} X-\nu_i\right)})}+m_i \ln \left(1+\|X\|^2\right),
  \end{equation}
where $h_{i}$, $\nu_{i}$, $m_{i}\in\mathbb{R}$, and each element of $\xi_{i}\in\mathbb{R}^{50}$ is randomly generated from the standard Gaussian distribution, for each $i\in\mathcal{V}$.

\lx{To the best of our knowledge, compressed distributed nonconvex optimization that can adapt to compressors with bounded relative compression error is only found in \cite{yi2022TAC, zhao2022NIPS}, and for compressors with globally and locally bounded absolute compression errors, it is only found in \cite{yi2022TAC}.
Hence, we compare Algorithms~1 and 2 with the BEER algorithm in \cite{zhao2022NIPS}, each paired with a compressor that has a bounded relative compression error.
Moreover, we compare Algorithms~1--3 with the compressed primal--dual algorithms (Algorithms~1--3 in \cite{yi2022TAC}), using three types of compressors: relative compression error, globally absolute error, and locally absolute error.}

(i) For the compressor with relative compression errors, we choose the classical norm-sign compressor as presented in \cite{Liao2022TAC}, and its mathematical model is given as follows:
\begin{align*}
\mathcal{C}(X)=\frac{\|X\|_{\infty}}{2} \operatorname{sign}(X),
\end{align*}
where this compressor is biased and non-contractive, yet it satisfies Assumption~\ref{assum-1} with $r=d/2$ and $\psi=1/d^{2}$, as demonstrated in \cite{Liao2022TAC}.
The algorithm parameters for comparison are set as follows.\\
%
$\bullet$ Algorithm~\ref{alg:1} with $\eta=0.8$, $\gamma=0.2$, $\varphi_{X}=0.3$, and $\varphi_{Y}=0.1$.\\
$\bullet$ Algorithm~\ref{alg:2} with $\eta=0.8$, $\gamma=0.2$, $\varphi_{X}=0.3$, $\varphi_{Y}=0.1$, and $\varsigma=0.3$.\\
$\bullet$ Algorithm~1 in \cite{yi2022TAC} with $\eta=0.3$, $\alpha=0.1$, $\beta=0.1$, and $\varphi=0.3$.\\
$\bullet$ Algorithm~2 in \cite{yi2022TAC} with $\eta=0.3$, $\alpha=0.1$, $\beta=0.1$, $\varphi=0.3$, and $\sigma=0.2$.\\
\XL{$\bullet$ BEER in \cite{zhao2022NIPS} with $\eta=1$, $\gamma=0.5$.}\\

\begin{figure}[!hbt]
\centering
\includegraphics[scale=0.8]{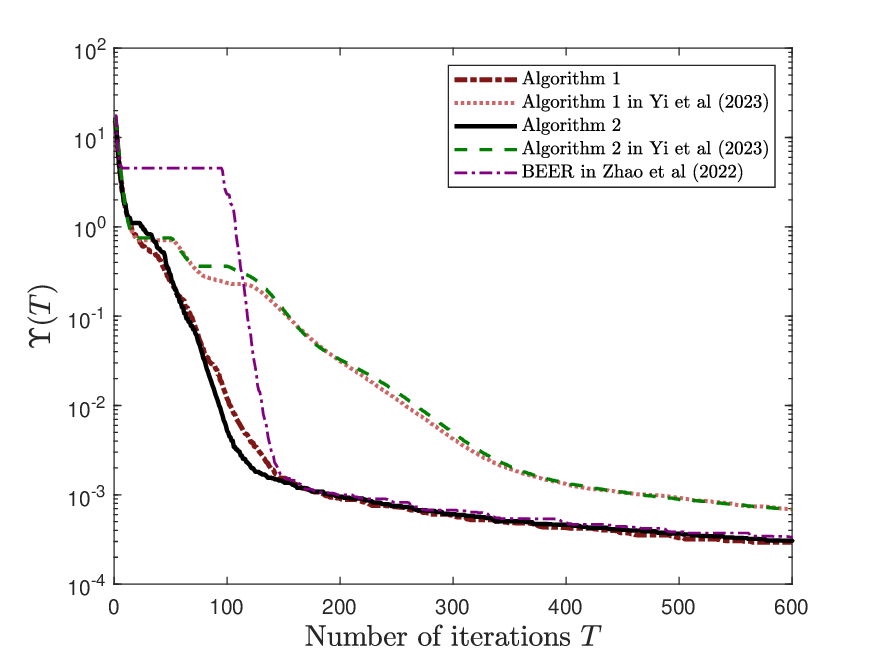}
\caption{Evolutions of $\Upsilon(T)$ with respect to the number of iterations of different compressed distributed algorithms with the norm-sign compressor.}
\label{Algorithm12C1}
\end{figure}

Fig.~\ref{Algorithm12C1} illustrates the evolution of
$$
\Upsilon(T)=\min_{k\in[T]}\{\sum_{i=1}^{n}\|X_{i}(k)-\bar{X}(k)\|^{2}+n\|\nabla F(\bar{X}(k))\|^{2}\},
$$
with respect to the number of iterations $k$, for the various compressed distributed algorithms.
\lx{From Fig. \ref{Algorithm12C1}, when utilizing the same norm-sign compressor, it is evident that the proposed Algorithms \ref{alg:1} and \ref{alg:2} demonstrate faster convergence speeds compared to Algorithms~1 and 2 in \cite{yi2022TAC}.}
Moreover, Algorithm~\ref{alg:2} exhibits a faster convergence speed than Algorithm~\ref{alg:1}, as it utilizes error feedback to correct the bias induced by biased compressors.

(ii) For the compressor with globally absolute compression errors, we choose the \lx{standard uniform quantizer}, and its mathematical model is as follows:
$$
\mathcal{C}(X)=\Delta\lfloor\frac{X}{\Delta}+\frac{\mathbf{1}_d}{2}\rfloor,
$$
where this compressor satisfies Assumption~\ref{assum-2} with $p=\infty$ and $C=0.25\Delta^{2}$, as demonstrated in \cite{yi2022TAC}.
In this section, we select $\Delta=2$.
The algorithm parameters for comparison are set as follows.\\
%
$\bullet$ Algorithm~\ref{alg:3} with $\eta=0.4$, and $\gamma=0.6$.\\
$\bullet$ Algorithm~\ref{alg:3} in \cite{yi2022TAC} with $\eta=0.25$, $\alpha=0.1$, and $\beta=0.1$.

\begin{figure}[!hbt]
\centering
\includegraphics[scale=0.8]{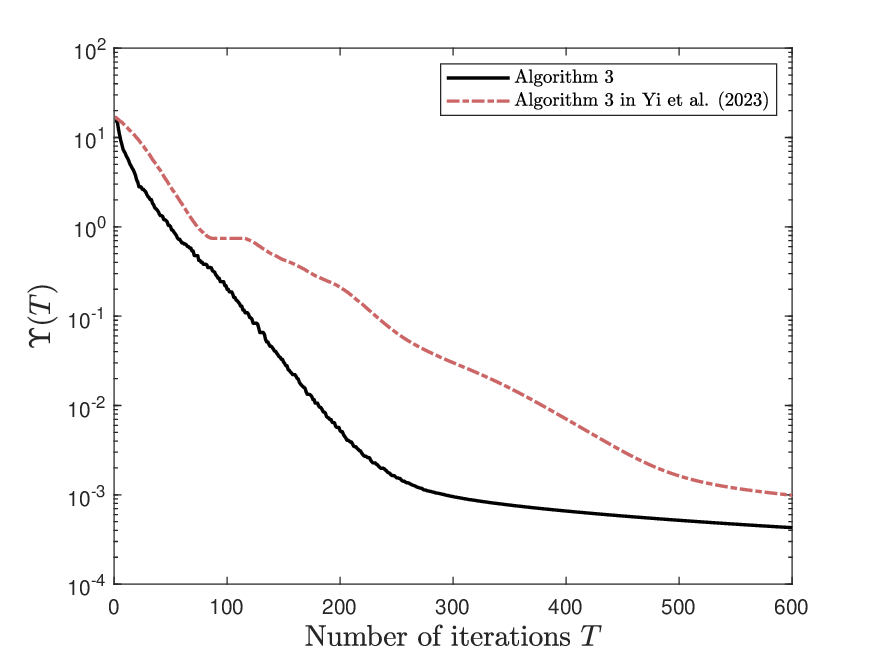}
\caption{Evolutions of $\Upsilon(T)$ with respect to the number of iterations of different compressed distributed algorithms with the standard uniform quantizer.}
\label{Algorithm3C2}
\end{figure}

Fig.~\ref{Algorithm3C2} depicts the evolution of $\Upsilon(T)$ in relation to the number of iterations $k$ for the different compressed distributed algorithms.
From Fig.\ref{Algorithm3C2}, when utilizing the same standard uniform quantizer, it is evident that the proposed Algorithm~\ref{alg:3} presents faster convergence compared to Algorithm~\ref{alg:3} in \cite{yi2022TAC}.

(iii) For the compressor with the locally absolute compression errors. We select the classical compressor: 1-bit binary quantizer, and the mathematical model is as follows:
$$
\mathcal{C}[X]=(q[X_{1}],\ldots,q[X_{d}]),
$$
and
  \begin{align*}
  q[a]&=
  \left\{\begin{aligned}
  &0.5,~a\geq0,\\
  &-0.5,~\text{otherwise},
  \end{aligned}\right.
  \end{align*}
where this compressor satisfies Assumption~\ref{assum-3} with $p=\infty$ and $\psi=0.5$.
The algorithm parameters for comparison are set as follows.\\
$\bullet$ Algorithm~\ref{alg:3} with $\eta=0.4$, and $\gamma=0.6$.\\
$\bullet$ Algorithm~\ref{alg:3} in \cite{yi2022TAC} with $\eta=0.25$, $\alpha=0.1$, and $\beta=0.1$.

\begin{figure}[!hbt]
\centering
\includegraphics[scale=0.8]{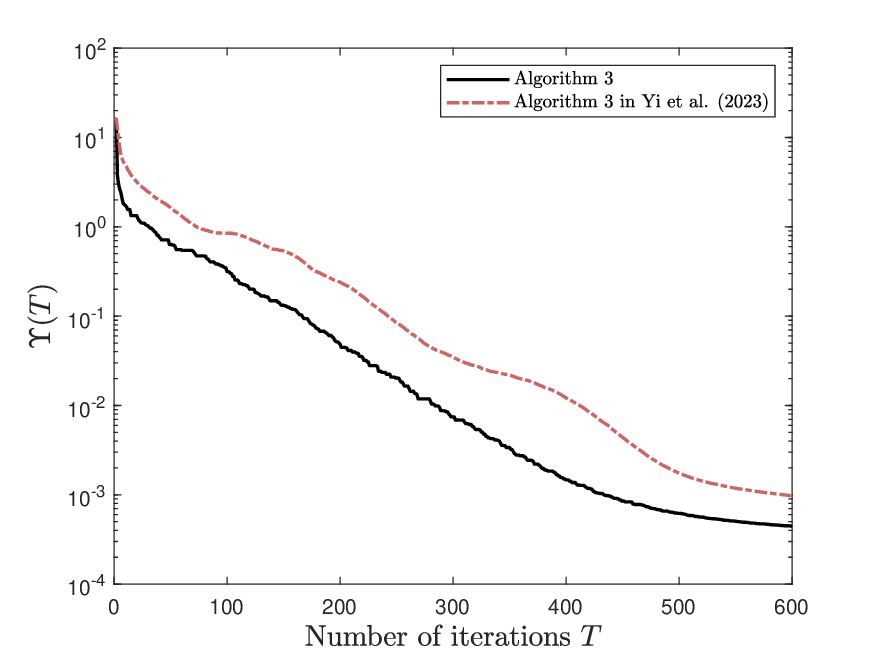}
\caption{Evolutions of $\Upsilon(T)$ with respect to the number of iterations of different compressed distributed algorithms with the 1-bit binary quantizer.}
\label{Algorithm3C3}
\end{figure}

Fig.~\ref{Algorithm3C3} illustrates the evolution of $\Upsilon(T)$ with respect to the number of iterations $k$ for the different compressed distributed algorithms.
From Fig.\ref{Algorithm3C3}, when utilizing the same 1-bit quantizer, it is evident that the proposed Algorithm~\ref{alg:3} demonstrates fast convergence speeds compared to Algorithm~3 in \cite{yi2022TAC}.


To compare the transmitted bits and $\mathcal{P}$ (percentage of the number of bits transmitted by the DGT algorithm proposed in \cite{Nedic2017SIAM}) for various algorithms and compressor combinations, achieving $\Upsilon(k)\leq10^{-3}$, Table~\ref{tbl:table4} is provided. This table presents a comparison of transmitted bits and $\mathcal{P}$ for different algorithms and compressor combinations, ensuring that $\Upsilon(k)\leq10^{-3}$. The experiment employs the proposed algorithm parameters outlined in (i)--(iii) and the DGT algorithm with $\eta=0.8$ and $\gamma=0.3$.
For the norm-sign compressor, transmitting one variable requires $2d+b_{C}$ bits if a scalar can be transmitted with $b_{C}$ bits with adequate precision. Here, we choose $b_{C}=64$. In the case of the standard uniform quantizer, transmitting one variable demands $db_{Q}$ bits if $b_{Q}$ bits are allocated to transmit an integer; here, we set $b_{Q}=4$. As for the 1-bit binary quantizer, transmitting one variable necessitates $d$ bits.
Table~\ref{tbl:table4} demonstrates that our algorithms require significantly fewer bits compared to the DGT algorithm to achieve a specific error level. Notably, the 1-bit binary quantizer uses only 4.35\% of the bits required by the DGT algorithm to achieve the same level of error.
Moreover, Algorithm~3 in \cite{yi2022TAC} requires only 2.94\% of the bits needed by the DGT algorithm to achieve $\Upsilon(T) \leq 10^{-3}$, making it the algorithm with the lowest number of transmitted bits.
\lx{Compared with the compressed distributed primal--dual algorithm in \cite{yi2022TAC}, in this paper, the proposed algorithms require each agent $i$ to communicate one additional $d$-dimensional variable $Q_{i}^{Y}(k)$ at each iteration, in addition to communicating $Q_{i}^{X}(k)$ with its neighbors.
The proposed algorithm is faster than Algorithm~3 in \cite{yi2022TAC}, but the latter has the lowest number of transmitted bits in Table~1.}

\begin{table}[!hbt]
\caption{Transmitted bits for different algorithms and compressors to reach $\Upsilon(T)\leq10^{-3}$.}
\centering
\resizebox{13cm}{!}
{
\begin{tabular}{|c|c|c|c|}
\hline
Algorithm & Compressor & Transmitted Bits & $\mathcal{P}$\\
\hline
Algorithm 1 & Norm--sign compressor & 61992 & 6.25\% \\
\hline
Algorithm 2 & Norm-sign compressor & 131856 & 13.29\% \\
\hline
Algorithm 3 & Standard uniform quantizer & 115600 & 11.65\% \\
\hline
Algorithm 3 & 1-bit binary quantizer & 43200 & 4.35\% \\
\hline
Algorithm 1 in \cite{yi2022TAC} & Norm--sign compressor & 74456 & 7.51\% \\
\hline
Algorithm 2 in \cite{yi2022TAC} & Norm-sign compressor & 152848 & 15.41\% \\
\hline
Algorithm 3 in \cite{yi2022TAC} & Standard uniform quantizer & 118200 & 11.92\% \\
\hline
Algorithm 3 in \cite{yi2022TAC} & 1-bit binary quantizer & 29200 & 2.94\% \\
\hline
\XL{BEER in \cite{zhao2022NIPS}}& \XL{Norm--sign compressor} & \XL{139072} & \XL{14.02\%} \\
\hline
\TY{DGT in \cite{Nedic2017SIAM}} & - & 992000 & - \\
\hline
\end{tabular}}
\label{tbl:table4}
\end{table}

\section{Conclusions}\label{sec-6}
In this paper, we introduced three classes of compressors to reduce the communication overhead. By integrating them with
DGT algorithm, we then proposed three distributed algorithms with compressed communication for distributed nonconvex optimization.
For the case where local cost functions are smooth, we designed several Lyapunov functions to demonstrate that the proposed algorithms sublinearly converge to a stationary point.
Moreover, when the global cost function satisfies the P--{\L} condition, we demonstrated that the proposed algorithm converges linearly to a global optimal point.
One future direction is to investigate \blue{general unbalanced directed graphs}.
\bibliographystyle{IEEEtran}
\bibliography{Ref}

\appendix
\hspace{-3mm}\emph{A. Useful Lemmas}\label{Appendix-A}

The following lemmas are used in the proofs.

Denote $\bm{B}=\mathrm{col}(B_{1},\ldots,B_{n})$, $\bm{D}=\mathrm{col}(D_{1},\ldots,D_{n})$, and $\bm{W}=W\otimes\bm{\mathrm{I}}_{d}$. Then, the compact form of Algorithm~\ref{alg:1} is
  \begin{subequations}\label{Alg1CMP}
  \begin{align}
  \bm{A}(k+1)&=\bm{A}(k)+\varphi_{X}\bm{Q}^{X}(k),\label{Alg1CMPb}\\
  \bm{B}(k+1)&=\bm{B}(k)+\varphi_{X}(\bm{Q}^{X}(k)-\bm{W}\bm{Q}^{X}(k)),\label{Alg1CMPc}\\
  \bm{C}(k+1)&=\bm{C}(k)+\varphi_{Y}\bm{Q}^{Y}(k),\label{Alg1CMPe}\\
  \bm{D}(k+1)&=\bm{D}(k)+\varphi_{Y}(\bm{Q}^{Y}(k)-\bm{W}\bm{Q}^{Y}(k)),\label{Alg1CMPf}\\
  \bm{X}(k+1)&=\bm{X}(k)-\gamma(\bm{B}(k)+(\bm{\mathrm{I}}_{nd}-\bm{W})\bm{Q}^{X}(k))-\eta\bm{Y}(k),\label{Alg1CMPg}\\
  \bm{Y}(k+1)&=\bm{Y}(k)-\gamma(\bm{D}(k)+(\bm{\mathrm{I}}_{nd}-\bm{W})\bm{Q}^{Y}(k))\notag\\
  &\quad+\nabla\tilde{F}(\bm{X}(k+1))-\nabla\tilde{F}(\bm{X}(k)),\label{Alg1CMPh}\\
  \bm{Q}^{X}(k+1)&=\mathcal{C}(\bm{X}(k+1)-\bm{A}(k+1)),\label{Alg1CMPa}\\
  \bm{Q}^{Y}(k+1)&=\mathcal{C}(\bm{Y}(k+1)-\bm{C}(k+1)),\label{Alg1CMPd}
  \end{align}
  \end{subequations}
where $\nabla\tilde{F}(\bm{X}(k))=[\nabla F^{T}_{1}(X_{1}(k)),\dots,\nabla F^{T}_{n}(X_{n}(k))]^{T}$.

For all $i\in[n]$, noted that $A_{i}(0)=B_{i}(0)=C_{i}(0)=D_{i}(0)=\bm{0}_{d}$, by mathematical induction, it is straightforward to check that $\bm{B}(k)=(\bm{\mathrm{I}}_{nd}-\bm{W})\bm{A}(k)$, and $\bm{D}(k)=(\bm{\mathrm{I}}_{nd}-\bm{W})\bm{C}(k)$.

We denote
  \begin{subequations}\label{Compression}
  \begin{align}
  \hat{\bm{X}}(k)&=\bm{A}(k)+\bm{Q}^{X}(k),\label{Compressiona}\\
  \hat{\bm{Y}}(k)&=\bm{C}(k)+\bm{Q}^{Y}(k),\label{Compressionb}
  \end{align}
  \end{subequations}
then, (\ref{Alg1CMPg}) and (\ref{Alg1CMPh}) respectively can be rewritten as
  \begin{subequations}\label{reAlg1}
  \begin{align}
  \bm{X}(k+1)&=\bm{X}(k)-\gamma(\bm{\mathrm{I}}_{nd}-\bm{W})\bm{\hat{X}}(k)-\eta\bm{Y}(k),\label{reAlg1a}\\
  \bm{Y}(k+1)&=\bm{Y}(k)-\gamma(\bm{\mathrm{I}}_{nd}-\bm{W})\bm{\hat{Y}}(k)+\nabla\tilde{F}(\bm{X}(k+1))-\nabla\tilde{F}(\bm{X}(k)).\label{reAlg1b}
  \end{align}
  \end{subequations}

\XL{If we replace $\hat{\bm{X}}(k)$ and $\hat{\bm{Y}}(k)$ with $\bm{X}(k)$ and $\bm{Y}(k)$ in (\ref{reAlg1}),}
then the DGT algorithm proposed in \cite{Nedic2017SIAM,Qu2017TCNS,Song2021Arxiv,Xiong2022TAC,Ma2021Springer} is derived.

\begin{lemma}\label{lem-spectral norm}
Suppose Assumption~\ref{assum-4} holds, then for all $\omega\in\mathbb{R}^{nd}$ and $\bar{\omega}=\bm{H}\omega$, we have
\begin{align}\label{lemEq}
&\|(\bm{\mathrm{I}}_{nd}+\gamma(\bm{W}-\bm{\mathrm{I}}_{nd}))(\omega-\bar{\omega})\|\notag\\
&=\|(\bm{\mathrm{I}}_{nd}+\gamma(\bm{W}-\bm{\mathrm{I}}_{nd})-\gamma\bm{H})(\omega-\bar{\omega})\|\notag\\
&\leq\|(\bm{\mathrm{I}}_{nd}-\gamma\bm{\mathrm{I}}_{nd})(\omega-\bar{\omega})\|+\|\gamma(\bm{W}-\bm{H})(\omega-\bar{\omega})\|\notag\\
&\leq(1-\gamma)\|\omega-\bar{\omega}\|+\sigma\gamma\|\omega-\bar{\omega}\|=\delta\|\omega-\bar{\omega}\|,
\end{align}
where $\delta=1-\gamma+\sigma\gamma$; the first equality holds due to $\bm{H}\omega=\bm{H}\bar{\omega}$; and the second inequality holds due to $\gamma\in(0,1)$.
\end{lemma}



\begin{lemma}\label{lem-barxy}
The following equalities hold:
\begin{align}\label{barx}
\bar{\bm{X}}(k+1)
&=\bm{H}\bm{X}(k+1)\notag\\
&=\bm{H}(\bm{X}(k)-\gamma(\bm{\mathrm{I}}_{nd}-\bm{W})\hat{\bm{X}}(k)-\eta\bm{Y}(k))\notag\\
&=\bar{\bm{X}}(k)-\eta\bar{\bm{Y}}(k),
\end{align}
where the second equality holds due to (\ref{reAlg1a}); and the last equality holds due to $\bm{H}(\bm{\mathrm{I}}_{nd}-\bm{W})=\bm{0}_{nd}$.
\begin{align}\label{bary}
\bar{\bm{Y}}(k+1)=&\bm{H}\bm{Y}(k+1)\notag\\
=&\bm{H}(\bm{Y}(k)-\gamma(\bm{\mathrm{I}}_{nd}-\bm{W})\hat{\bm{Y}}(k)+\bm{g}(k+1)-\bm{g}(k))\notag\\
=&\bar{\bm{g}}(k+1),
\end{align}
where $\bm{g}(k)=\nabla\tilde{F}(\bm{X}(k))$, $\tilde{F}(\bm{X})=\sum_{i=1}^{n}F_{i}(X_{i})$, and $\bar{\bm{g}}(k)=\bm{H}\bm{g}(k)$;
the second equality holds due to (\ref{reAlg1b}); the last equality holds due to $\bm{H}(\bm{\mathrm{I}}_{nd}-\bm{W})=\bm{0}_{nd}$ and $\bm{Y}(0)=\nabla\tilde{F}(\bm{X}(0))$.
\end{lemma}

\hspace{-3mm}\emph{B. Proof of Theorem~\ref{thm-1}}

For the $\|\bm{X}(k+1)-\bar{\bm{X}}(k+1)\|$, we have
\begin{align}\label{ErrorxbarxnoE}
&\|\bm{X}(k+1)-\bar{\bm{X}}(k+1)\|^{2}\notag\\
&=\|\bm{X}(k)-\gamma(\bm{\mathrm{I}}_{nd}-\bm{W})\bm{\hat{X}}(k)-\eta\bm{Y}(k)-(\bar{\bm{X}}(k)-\eta\bar{\bm{Y}}(k))\|^{2}\notag\\
&=\|\gamma(\bm{\mathrm{I}}_{nd}-\bm{W})(\bm{X}(k)-\bm{\hat{X}}(k))+(\bm{\mathrm{I}}_{nd}+\gamma(\bm{W}-\bm{\mathrm{I}}_{nd}))\bm{X}(k)-\bar{\bm{X}}(k)-\eta(\bm{Y}(k)-\bar{\bm{Y}}(k))\|^{2}\notag\\
&=\|\gamma(\bm{\mathrm{I}}_{nd}-\bm{W})(\bm{X}(k)-\bm{\hat{X}}(k))+(\bm{\mathrm{I}}_{nd}+\gamma(\bm{W}-\bm{\mathrm{I}}_{nd}))(\bm{X}(k)-\bar{\bm{X}}(k))-\eta(\bm{Y}(k)-\bar{\bm{Y}}(k))\|^{2}\notag\\
&\leq(2\gamma\|\bm{X}(k)-\bm{\hat{X}}(k)\|+\delta\|\bm{X}(k)-\bar{\bm{X}}(k)\|+\eta\|\bm{Y}(k)-\bar{\bm{Y}}(k)\|)^{2}\notag\\
&\leq\frac{8\gamma}{1-\sigma}\|\bm{X}(k)-\bm{\hat{X}}(k)\|^{2}+\delta\|\bm{X}(k)-\bar{\bm{X}}(k)\|^{2}+\frac{2\eta^{2}}{\gamma(1-\sigma)}\|\bm{Y}(k)-\bar{\bm{Y}}(k)\|^{2},
\end{align}
where the first equality holds due to (\ref{reAlg1a}) and (\ref{barx}); the second equality holds due to (\ref{bary}); the last equality holds due to $(\bm{W}-\bm{\mathrm{I}}_{nd})\bm{H}=\bm{0}_{nd}$; the first inequality holds due to (\ref{lemEq}); and the last inequality holds due to the Cauchy--Schwarz inequality.

Based on (\ref{Assum1a}), (\ref{Alg1CMPa}) and (\ref{Compressiona}), we have
\begin{equation}\label{ExpErrorx}
\bm{\mathrm{E}}_{\mathcal{C}}[\|\bm{X}(k)-\hat{\bm{X}}(k)\|^{2}]
=\bm{\mathrm{E}}_{\mathcal{C}}[\|\bm{X}(k)-\bm{A}(k)-\mathcal{C}(\bm{X}(k)-\bm{A}(k))\|^{2}]
\leq C\|\bm{X}(k)-\bm{A}(k)\|^{2}.
\end{equation}

From (\ref{ErrorxbarxnoE}) and (\ref{ExpErrorx}), we have
\begin{align}\label{Errorxbarx}
&\bm{\mathrm{E}}_{\mathcal{C}}[\|\bm{X}(k+1)-\bar{\bm{X}}(k+1)\|^{2}]\notag\\
&\leq\bm{\mathrm{E}}_{\mathcal{C}}[\delta\|\bm{X}(k)-\bar{\bm{X}}(k)\|^{2}+\frac{2\eta^{2}}{\gamma(1-\sigma)}\|\bm{Y}(k)-\bar{\bm{Y}}(k)\|^{2}+\frac{8\gamma}{1-\sigma}C\|\bm{X}(k)-\bm{A}(k)\|^{2}],
\end{align}
where the inequality holds due to the Cauchy--Schwarz inequality.



For the $\|\bm{Y}(k+1)-\bar{\bm{Y}}(k+1)\|$, we have
\begin{align}\label{ErrorygnoE}
&\|\bm{Y}(k+1)-\bar{\bm{Y}}(k+1)\|^{2}\notag\\
&=\|\bm{Y}(k)-\gamma(\bm{\mathrm{I}}_{nd}-\bm{W})\bm{\hat{Y}}(k)+\bm{g}(k+1)-\bm{g}(k)-\bar{\bm{g}}(k+1)\|^{2}\notag\\
&=\|\gamma(\bm{\mathrm{I}}_{nd}-\bm{W})(\bm{Y}(k)-\bm{\hat{Y}}(k))+(\bm{\mathrm{I}}_{nd}+\gamma(\bm{W}-\bm{\mathrm{I}}_{nd}))\notag\\
&\quad\times(\bm{Y}(k)-\bar{\bm{Y}}(k))+\bm{g}(k+1)-\bm{g}(k)-\bar{\bm{g}}(k+1)+\bar{\bm{g}}(k)\|^{2}\notag\\
&=\|\gamma(\bm{\mathrm{I}}_{nd}-\bm{W})(\bm{Y}(k)-\bm{\hat{Y}}(k))+(\bm{\mathrm{I}}_{nd}+\gamma(\bm{W}-\bm{\mathrm{I}}_{nd}))\notag\\
&\quad\times(\bm{Y}(k)-\bar{\bm{Y}}(k))+(\bm{\mathrm{I}}_{nd}-\bm{H})(\bm{g}(k+1)-\bm{g}(k))\|^{2}\notag\\
&\leq(2\gamma\|\bm{Y}(k)-\bm{\hat{Y}}(k)\|+\delta\|\bm{Y}(k)-\bar{\bm{Y}}(k)\|+\bm{L}_{f}\|\bm{X}(k+1)-\bm{X}(k)\|)^{2}\notag\\
&\leq\frac{8\gamma}{1-\sigma}\|\bm{Y}(k)-\bm{\hat{Y}}(k)\|^{2}+\delta\|\bm{Y}(k)-\bar{\bm{Y}}(k)\|^{2}+\frac{2 L_{f}^{2}}{\gamma(1-\sigma)}\|\bm{X}(k+1)-\bm{X}(k)\|^{2},
\end{align}
where the first equality holds due to (\ref{reAlg1b}); the second equality holds due to $(\bm{W}-\bm{\mathrm{I}}_{nd})\bm{H}=\bm{0}_{nd}$; the last equality holds due to (\ref{bary});
the first inequality holds due to (\ref{assum2}), (\ref{lemEq}), and $\rho(\bm{\mathrm{I}}_{nd}-\bm{H})\leq1$; the last inequality holds due to (\ref{bary}) and the Cauchy--Schwarz inequality.

Denote $\bm{g}^{0}(k)=\nabla\tilde{F}(\bar{\bm{X}}(k))$, and $\bar{\bm{g}}^{0}(k)=\bm{H}\bm{g}^{0}(k)=\bm{1}_{n}\otimes\nabla F(\bar{X}(k))$.
Then, from Assumption~\ref{assum-6} and $\rho(\bm{H})=1$, we have
\begin{equation}\label{barg}
\|\bar{\bm{g}}^{0}(k)-\bar{\bm{g}}(k)\|=\|\bm{H}(\bm{g}^{0}(k)-\bm{g}(k))\leq\|\bm{g}^{0}(k)-\bm{g}(k)\|\leq L_{f}\|\bm{X}(k)-\bar{\bm{X}}(k)\|.
\end{equation}

Based on (\ref{reAlg1a}), (\ref{bary}), and (\ref{barg}), it can be obtained that
\begin{align}\label{Errorxx}
&\|\bm{X}(k+1)-\bm{X}(k)\|\notag\\
&=\|\bm{X}(k)-\gamma(\bm{\mathrm{I}}_{nd}-\bm{W})\bm{\hat{X}}(k)-\eta\bm{Y}(k)-\bm{X}(k)\|\notag\\
&=\|\gamma(\bm{\mathrm{I}}_{nd}-\bm{W})(\bm{X}(k)-\bm{\hat{X}}(k))+\gamma(\bm{W}-\bm{\mathrm{I}}_{nd})(\bm{X}(k)-\bar{\bm{X}}(k))-\eta(\bm{Y}(k)-\bar{\bm{Y}}(k))\notag\\
&\quad+\eta(\bar{\bm{g}}^{0}(k)-\bar{\bm{g}}(k))-\eta\bar{\bm{g}}^{0}(k)\|\notag\\
&\leq2\gamma\|\bm{X}(k)-\bm{\hat{X}}(k)\|+(2\gamma+\eta L_{f})\|\bm{X}(k)-\bar{\bm{X}}(k)\|+\eta\|\bm{Y}(k)-\bar{\bm{Y}}(k)\|+\eta\|\bar{\bm{g}}^{0}(k)\|,
\end{align}
where the second equality holds due to $(\bm{W}-\bm{\mathrm{I}}_{nd})\bm{H}=\bm{0}_{nd}$.

Based on (\ref{Assum1a}), (\ref{Alg1CMPd}) and (\ref{Compressionb}), it can be calculated that
\begin{equation}\label{ExpErrory}
\bm{\mathrm{E}}_{\mathcal{C}}[\|\bm{Y}(k)-\hat{\bm{Y}}(k)\|^{2}]=\bm{\mathrm{E}}_{\mathcal{C}}[\|\bm{Y}(k)-\bm{C}(k)-\mathcal{C}(\bm{Y}(k)-\bm{C}(k))\|^{2}]\leq C\|\bm{Y}(k)-\bm{C}(k)\|^{2}.
\end{equation}

From (\ref{ExpErrorx}), (\ref{ErrorygnoE}), (\ref{Errorxx}), and (\ref{ExpErrory}), we have
\begin{align}\label{Erroryg}
&\bm{\mathrm{E}}_{\mathcal{C}}[\|\bm{Y}(k+1)-\bar{\bm{Y}}(k+1)\|^{2}]\notag\\
&\leq\bm{\mathrm{E}}_{\mathcal{C}}[\frac{8\gamma}{1-\sigma}\|\bm{Y}(k)-\bm{\hat{Y}}(k)\|^{2}+\delta\|\bm{Y}(k)-\bar{\bm{Y}}(k)\|^{2}+\frac{8L_{f}^{2}}{\gamma(1-\sigma)}(4\gamma^{2}\|\bm{X}(k)-\bm{\hat{X}}(k)\|^{2}\notag\\
&\quad+(8\gamma^{2}+2\eta^{2} L_{f}^{2})\|\bm{X}(k)-\bar{\bm{X}}(k)\|^{2}+\eta^{2}\|\bm{Y}(k)-\bar{\bm{Y}}(k)\|^{2}+\eta^{2}\|\bar{\bm{g}}^{0}(k)\|^{2})]\notag\\
&\leq\bm{\mathrm{E}}_{\mathcal{C}}[(\delta+\frac{8L_{f}^{2}}{\gamma(1-\sigma)}\eta^{2})\|\bm{Y}(k)-\bar{\bm{Y}}(k)\|^{2}+\frac{8L_{f}^{2}}{\gamma(1-\sigma)}(4\gamma^{2}C\|\bm{X}(k)-\bm{A}(k)\|^{2}\notag\\
&\quad+(8\gamma^{2}+2\eta^{2} L_{f}^{2})\|\bm{X}(k)-\bar{\bm{X}}(k)\|^{2}+\eta^{2}\|\bar{\bm{g}}^{0}(k)\|^{2})+\frac{8\gamma}{1-\sigma}C\|\bm{Y}(k)-\bm{C}(k)\|^{2}].
\end{align}

For the $\bm{\mathrm{E}}_{\mathcal{C}}[\|\bm{X}(k+1)-\bm{A}(k+1)\|^{2}]$, we have
\begin{align}\label{ExpErrorxa1}
&\bm{\mathrm{E}}_{\mathcal{C}}[\|\bm{X}(k+1)-\bm{A}(k+1)\|^{2}]\notag\\
&=\bm{\mathrm{E}}_{\mathcal{C}}[\|\bm{X}(k+1)-\bm{X}(k)+\bm{X}(k)-\bm{A}(k)-\varphi_{X}\bm{Q}^{X}(k)\|^{2}]\notag\\
&=\bm{\mathrm{E}}_{\mathcal{C}}[\|\bm{X}(k+1)-\bm{X}(k)+(1-\varphi_{X}r)(\bm{X}(k)-\bm{A}(k))\notag\\
&\quad+\varphi_{X}r((\bm{X}(k)-\bm{A}(k))-\mathcal{C}_{r}(\bm{X}(k)-\bm{A}(k)))\|^{2}]\notag\\
&\leq\bm{\mathrm{E}}_{\mathcal{C}}[(1+c_{1}^{-1})\|\bm{X}(k+1)-\bm{X}(k)\|^{2}+(1+c_{1})(1-\varphi_{X}r)\|\bm{X}(k)-\bm{A}(k)\|^{2}\notag\\
&\quad+(1+c_{1})\varphi_{X}r\|\bm{X}(k)-\bm{A}(k)-\mathcal{C}_{r}(\bm{X}(k)-\bm{A}(k))\|^{2}]\notag\\
&\leq\bm{\mathrm{E}}_{\mathcal{C}}[(1+c_{1}^{-1})\|\bm{X}(k+1)-\bm{X}(k)\|^{2}+(1+c_{1})(1-\varphi_{X}r)\|\bm{X}(k)-\bm{A}(k)\|^{2}\notag\\
&\quad+(1+c_{1})\varphi_{X}r(1-\psi)\|\bm{X}(k)-\bm{A}(k)\|^{2}]\notag\\
&=\bm{\mathrm{E}}_{\mathcal{C}}[(1+c_{1}^{-1})\|\bm{X}(k+1)-\bm{X}(k)\|^{2}+(1-c_{1}-2c_{1}^{2})\|\bm{X}(k)-\bm{A}(k)\|^{2}]\notag\\
&\leq\bm{\mathrm{E}}_{\mathcal{C}}[(1+c_{1}^{-1})(2\gamma\|\bm{X}(k)-\bm{\hat{X}}(k)\|+(2\gamma+\eta L_{f})\|\bm{X}(k)-\bar{\bm{X}}(k)\|+\eta\|\bm{Y}(k)-\bar{\bm{Y}}(k)\|\notag\\
&\quad+\eta\|\bar{\bm{g}}^{0}(k)\|)^{2}+(1-c_{1}-2c_{1}^{2})\|\bm{X}(k)-\bm{A}(k)\|^{2}]\notag\\
&\leq\bm{\mathrm{E}}_{\mathcal{C}}[4(1+c_{1}^{-1})(4\gamma^{2}C\|\bm{X}(k)-\bm{A}(k)\|^{2}+(8\gamma^{2}+2\eta^{2}L_{f}^{2})\|\bm{X}(k)-\bar{\bm{X}}(k)\|^{2}\notag\\
&\quad+\eta^{2}\|\bm{Y}(k)-\bar{\bm{Y}}(k)\|^{2}+\eta^{2}\|\bar{\bm{g}}^{0}(k)\|^{2})+(1-c_{1}-2c_{1}^{2})\|\bm{X}(k)-\bm{A}(k)\|^{2}],
\end{align}
where $\mathcal{C}_{r}(\cdot)=\frac{\mathcal{C}(\cdot)}{r}$;
the first and second equalities hold due to (\ref{Alg1CMP});
the first inequality holds due to the Cauchy--Schwarz inequality;
the second inequality holds due to (\ref{Assum1b});
the third inequality holds due to (\ref{Errorxx});
the last inequality holds due to (\ref{ExpErrorx}), and the Cauchy--Schwarz inequality.


For the $\bm{\mathrm{E}}_{\mathcal{C}}[\|\bm{Y}(k+1)-\bm{C}(k+1)\|^{2}]$, we have
\begin{align}\label{ExpErroryc2a}
&\bm{\mathrm{E}}_{\mathcal{C}}[\|\bm{Y}(k+1)-\bm{C}(k+1)\|^{2}]\notag\\
&=\bm{\mathrm{E}}_{\mathcal{C}}[\|\bm{Y}(k+1)-\bm{Y}(k)+\bm{Y}(k)-\bm{C}(k)-\varphi_{Y}\bm{Q}^{Y}(k)\|^{2}]\notag\\
&=\bm{\mathrm{E}}_{\mathcal{C}}[\|\bm{Y}(k+1)-\bm{Y}(k)+(1-\varphi_{Y}r)(\bm{Y}(k)-\bm{C}(k))\notag\\
&\quad+\varphi_{Y}r((\bm{Y}(k)-\bm{C}(k))-\mathcal{C}_{r}(\bm{Y}(k)-\bm{C}(k)))\|^{2}]\notag\\
&\leq\bm{\mathrm{E}}_{\mathcal{C}}[(1+c_{2}^{-1})\|\bm{Y}(k+1)-\bm{Y}(k)\|^{2}+(1+c_{2})(1-\varphi_{Y}r)\|\bm{Y}(k)-\bm{C}(k)\|^{2}\notag\\
&\quad+(1+c_{2})\varphi_{Y}r\|\bm{Y}(k)-\bm{C}(k)-\mathcal{C}_{r}(\bm{Y}(k)-\bm{C}(k))\|^{2}]\notag\\
&\leq\bm{\mathrm{E}}_{\mathcal{C}}[(1+c_{2}^{-1})\|\bm{Y}(k+1)-\bm{Y}(k)\|^{2}+(1+c_{2})(1-\varphi_{Y}r)\|\bm{Y}(k)-\bm{C}(k)\|^{2}\notag\\
&\quad+(1+c_{2})\varphi_{Y}r(1-\psi)\|\bm{Y}(k)-\bm{C}(k)\|^{2}]\notag\\
&=\bm{\mathrm{E}}_{\mathcal{C}}[(1+c_{2}^{-1})\|\bm{Y}(k+1)-\bm{Y}(k)\|^{2}+(1-c_{2}-2c_{2}^{2})\|\bm{Y}(k)-\bm{C}(k)\|^{2}],
\end{align}
where the first equality holds due to (\ref{Alg1CMPe});
the second equality holds due to (\ref{Alg1CMPd});
the first inequality holds due to the Cauchy--Schwarz inequality;
and the last inequality holds due to (\ref{Assum1b}).

Based on (\ref{reAlg1b}), (\ref{barg}) and (\ref{Errorxx}), it can be derived that
\begin{align}\label{Erroryy}
&\|\bm{Y}(k+1)-\bm{Y}(k)\|\notag\\
&=\|\bm{Y}(k)-\gamma(\bm{\mathrm{I}}_{nd}-\bm{W})\bm{\hat{Y}}(k)+\bm{g}(k+1)-\bm{g}(k)-\bm{Y}(k)\|\notag\\
&=\|\gamma(\bm{\mathrm{I}}_{nd}-\bm{W})(\bm{Y}(k)-\bm{\hat{Y}}(k))+\gamma(\bm{W}-\bm{\mathrm{I}}_{nd})(\bm{Y}(k)-\bar{\bm{Y}}(k))+\bm{g}(k+1)-\bm{g}(k)\|\notag\\
&\leq2\gamma\|\bm{Y}(k)-\bm{\hat{Y}}(k)\|+2\gamma\|\bm{Y}(k)-\bar{\bm{Y}}(k)\|+L_{f}\|\bm{X}(k+1)-\bm{X}(k)\|\notag\\
&\leq2\gamma\|\bm{Y}(k)-\bm{\hat{Y}}(k)\|+(2\gamma+\eta L_{f})\|\bm{Y}(k)-\bar{\bm{Y}}(k)\|+L_{f}(2\gamma\|\bm{X}(k)-\bm{\hat{X}}(k)\|\notag\\
&\quad+(2\gamma+\eta L_{f})\|\bm{X}(k)-\bar{\bm{X}}(k)\|+\eta\|\bar{\bm{g}}^{0}(k)\|),
\end{align}
where the last equality holds due to $(\bm{W}-\bm{\mathrm{I}}_{nd})\bm{H}=\bm{0}_{nd}$ and (\ref{bary}).

From (\ref{ExpErroryc2a}) and (\ref{Erroryy}), we have
\begin{align}\label{ExpErroryc2}
&\bm{\mathrm{E}}_{\mathcal{C}}[\|\bm{Y}(k+1)-\bm{C}(k+1)\|^{2}]\notag\\
&\leq\bm{\mathrm{E}}_{\mathcal{C}}[(1+c_{2}^{-1})(2\gamma\|\bm{Y}(k)-\bm{\hat{Y}}(k)\|+(2\gamma+\eta L_{f})\|\bm{Y}(k)-\bar{\bm{Y}}(k)\|\notag\\
&\quad+L_{f}(2\gamma\|\bm{X}(k)-\bm{\hat{X}}(k)\|+(2\gamma+\eta L_{f})\|\bm{X}(k)-\bar{\bm{X}}(k)\|+\eta\|\bar{\bm{g}}^{0}(k)\|))^{2}\notag\\
&\quad+(1-c_{2}-2c_{2}^{2})\|\bm{Y}(k)-\bm{C}(k)\|^{2}]\notag\\
&\leq\bm{\mathrm{E}}_{\mathcal{C}}[5(1+c_{2}^{-1})(4\gamma^{2}C\|\bm{Y}(k)-\bm{C}(k)\|^{2}+(8\gamma^{2}+2\eta^{2}L_{f}^{2})\|\bm{Y}(k)-\bar{\bm{Y}}(k)\|^{2}\notag\\
&\quad+L_{f}^{2}(4\gamma^{2}C\|\bm{X}(k)-\bm{A}(k)\|^{2}+(8\gamma^{2}+2\eta^{2}L_{f}^{2})\|\bm{X}(k)-\bar{\bm{X}}(k)\|^{2}+\eta^{2}\|\bar{\bm{g}}^{0}(k)\|^{2}))\notag\\
&\quad+(1-c_{2}-2c_{2}^{2})\|\bm{Y}(k)-\bm{C}(k)\|^{2}],
\end{align}
where the last inequality holds due to (\ref{ExpErrorx}), (\ref{ExpErrory}), and the Cauchy--Schwarz inequality.

From \cite[Lemma~1.2.3]{Nesterov2018Lectures}, we know that (\ref{assum2}) implies
\begin{equation}\label{Lipschitz2}
|F_{i}(Y)-F_{i}(X)-(Y-X)^{T}\nabla F_{i}(X)|\leq\frac{L_{f}}{2}\|Y-X\|^{2},~\forall X,~Y\in\mathbb{R}^{d}.
\end{equation}

For the $n(F(\bar{X}(k+1)-F^{\star}))$, based on (\ref{Lipschitz2}), we have
\begin{align}\label{OptimalF2}
&n(F(\bar{X}(k+1))-F^{\star})\notag\\
&= n(F(\bar{X}(k))- F^{\star}+ F(\bar{X}(k+1))- F(\bar{X}(k)))\notag\\
&\leq n(F(\bar{X}(k))- F^{\star})-\eta\bar{\bm{g}}^{T}(k)\bm{g}^{0}(k)+\frac{\eta^{2}L_{f}}{2}\|\bar{\bm{g}}(k)\|^{2}\notag\\
&= n(F(\bar{X}(k))- F^{\star})-\eta\bar{\bm{g}}^{T}(k)\bar{\bm{g}}^{0}(k)+\frac{\eta^{2}L_{f}}{2}\|\bar{\bm{g}}(k)\|^{2}\notag\\
&=n(F(\bar{X}(k))-F^{\star})-\frac{\eta}{2}\bar{\bm{g}}^{T}(k)(\bar{\bm{g}}(k)+\bar{\bm{g}}^{0}(k)-\bar{\bm{g}}(k))\notag\\
&\quad-\frac{\eta}{2}(\bar{\bm{g}}(k)+\bar{\bm{g}}^{0}(k)-\bar{\bm{g}}^{0}(k))^{T}\bar{\bm{g}}^{0}(k)+\frac{\eta^{2}L_{f}}{2}\|\bar{\bm{g}}(k)\|^{2}\notag\\
&\leq n(F(\bar{X}(k))-F^{\star})-\frac{\eta}{4}\|\bar{\bm{g}}(k)\|^{2}+\frac{\eta}{4}\|\bar{\bm{g}}^{0}(k)-\bar{\bm{g}}(k)\|^{2}\notag\\
&\quad-\frac{\eta}{4}\|\bar{\bm{g}}^{0}(k)\|^{2}+\frac{\eta}{4}\|\bar{\bm{g}}^{0}(k)-\bar{\bm{g}}(k)\|^{2}+\frac{\eta^{2} L_{f}}{2}\|\bar{\bm{g}}(k)\|^{2}\notag\\
&\leq n(F(\bar{X}(k))-F^{\star})-\frac{\eta}{4}(1-2 \eta L_{f})\|\bar{\bm{g}}(k)\|^{2}\notag\\
&\quad+\frac{\eta}{2}L_{f}^{2}\|\bm{X}(k)-\bar{\bm{X}}(k)\|^{2}-\frac{\eta}{4}\|\bar{\bm{g}}^{0}(k)\|^{2},
\end{align}
where the first inequality holds since that $F$ is smooth and (\ref{Lipschitz2}); the second equality holds due to $\bar{\bm{g}}^{T}(k)\bm{g}^{0}(k)=\bm{g}^{T}(k)\bm{H}\bm{g}^{0}(k)=\bm{g}^{T}(k)\bm{H}\bm{H}\bm{g}^{0}(k)=\bar{\bm{g}}^{T}(k)\bar{\bm{g}}^{0}(k)$; the second inequality holds due to the Cauchy--Schwarz inequality; and the last inequality holds due to (\ref{barg}).

From (\ref{ErrorxbarxnoE}), (\ref{ErrorygnoE}), (\ref{ExpErrorxa1}), (\ref{ExpErroryc2a}), we have
\begin{align}\label{Lyapunovkp1}
&\bm{\mathrm{E}}_{\mathcal{C}}[U(k+1)]\notag\\
&\leq\bm{\mathrm{E}}_{\mathcal{C}}[\theta_{5}\|\bm{X}(k)-\bar{\bm{X}}(k)\|^{2}+\theta_{6}\|\bm{Y}(k)-\bar{\bm{Y}}(k)\|^{2}+\theta_{7}\|\bm{X}(k)-\bm{A}(k)\|^{2}\notag\\
&\quad+\theta_{8}\|\bm{Y}(k)-\bm{C}(k)\|^{2}+n(F(\bar{X}(k))-F^{\star})-\theta_{2}\|\bar{\bm{g}}^{0}(k)\|^{2}-\theta_{9}\|\bar{\bm{g}}(k)\|^{2}],
\end{align}
where
\begin{align*}
\theta_{5}&=\delta+\phi\xi_{1}+\xi_{2}+\xi_{3}+\frac{\eta L_{f}^{2}}{2},~\theta_{6}=\xi_{4}+\phi\xi_{5}+\xi_{6}+\xi_{7},\\
\theta_{7}&=\xi_{8}+\phi\xi_{9}+\xi_{10}+\xi_{11},~\theta_{8}=\phi\xi_{8}+\xi_{12},~\theta_{9}=\frac{\eta}{4}(1-2\eta L_{f}),\\
\xi_{1}&=\frac{8L_{f}^{2}}{(1-\sigma)\gamma}(8\gamma^{2}+2\eta^{2}L_{f}^{2}),~\xi_{2}=4(8\gamma^{2}+2\eta^{2}L_{f}^{2})(1+c_{1}^{-1}),~\xi_{3}=\xi_{7}L_{f}^{2},\\
\xi_{4}&=\frac{2\eta^{2}}{\gamma(1-\sigma)},~\xi_{5}=\delta+\frac{8L_{f}^{2}}{(1-\sigma)\gamma}\eta^{2},~\xi_{6}=4(1+c_{1}^{-1})\eta^{2},\\
\xi_{7}&=5(1+c_{2}^{-1})(8\gamma^{2}+2\eta^{2}L_{f}^{2}),~\xi_{8}=\frac{8\gamma}{1-\sigma}C,\\
\xi_{9}&=\frac{32L_{f}^{2}}{1-\sigma}\gamma C,~\xi_{10}=16\gamma^{2}(1+c_{1}^{-1})C+(1-c_{1}-2c_{1}^{2}),\\
\xi_{11}&=5(1+c_{2}^{-1})4\gamma^{2}L_{f}^{2}C,~\xi_{12}=20\gamma^{2}(1+c_{2}^{-1})C+(1-c_{2}-2c_{2}^{2}).
\end{align*}

Next, we prove that $\theta_{5},~\theta_{6},~\theta_{7},~\theta_{8}\in(0,1)$, $\theta_{2},~\theta_{9}\geq0$.

(i) We first prove that the following inequality holds.
\begin{equation}\label{omega1}
0<\theta_{5}<1-0.37(1-\sigma)\gamma<1.
\end{equation}

From $\eta<\frac{(1-\sigma)^{2}\gamma}{40L_{f}}$, we have
\begin{equation}\label{omega1a}
\phi\xi_{1}<\frac{(1-\sigma)^{2}}{320L_{f}^{2}}\frac{64.01L_{f}^{2}\gamma^{2}}{(1-\sigma)\gamma}<0.21(1-\sigma)\gamma.
\end{equation}

From $\gamma<\frac{1-\sigma}{160(1+c_{1}^{-1})}$, and $\eta<\frac{(1-\sigma)^{2}\gamma}{40L_{f}}$, we have
\begin{equation}\label{omega1b}
\xi_{2}<\frac{32.005(1-\sigma)}{160(1+c_{1}^{-1})}\gamma(1+c_{1}^{-1})<  0.21(1-\sigma)\gamma.
\end{equation}

From $\gamma<\frac{1-\sigma}{40000(1+c_{2}^{-1})L_{f}^{2}}$, and $\eta<\frac{(1-\sigma)^{2}\gamma}{40L_{f}}$, we have
\begin{equation}\label{omega1c}
\xi_{3}<40.00625\gamma^{2}L_{f}^{2}<0.0011(1-\sigma)\gamma.
\end{equation}

From $\eta<\frac{0.4(1-\sigma)\gamma}{L_{f}^{2}}$, we have
\begin{equation}\label{omega1d}
\eta L_{f}^{2}/2<0.2(1-\sigma)\gamma.
\end{equation}

From $\gamma<\frac{1-\sigma}{160(1+c_{1}^{-1})}$, and (\ref{omega1a})--(\ref{omega1d}), we know that (\ref{omega1}) holds. This completes the proof.

(ii) We second prove that the following inequality holds.
\begin{equation}\label{omega2}
0<\theta_{6}<\phi(1-0.07(1-\sigma)\gamma)<\phi.
\end{equation}

From $\eta<\frac{(1-\sigma)^{2}\gamma}{40L_{f}}$, we have
\begin{equation}\label{omega2a}
\xi_{4}<\frac{2}{\gamma(1-\sigma)}\frac{(1-\sigma)^{4}\gamma^{2}}{1600L_{f}^{2}}=\frac{(1-\sigma)^{2}}{320L_{f}^{2}}\frac{2(1-\sigma)\gamma}{5}=0.4\phi(1-\sigma)\gamma.
\end{equation}

From $\eta<\frac{(1-\sigma)^{2}\gamma}{40L_{f}}$, we have
\begin{equation}\label{omega2b}
4\phi\frac{2L_{f}^{2}}{(1-\sigma)\gamma}\eta^{2}<\phi\frac{8L_{f}^{2}}{(1-\sigma)\gamma}\frac{(1-\sigma)^{4}\gamma^{2}}{1600L_{f}^{2}}<0.005\phi(1-\sigma)\gamma.
\end{equation}

From $\eta<\frac{(1-\sigma)^{2}}{80L_{f}}\sqrt{\frac{\gamma}{1+c_{1}^{-1}}}$, we have
\begin{equation}\label{omega2c}
\xi_{6}<4\frac{(1-\sigma)^{4}}{6400L_{f}^{2}}\frac{\gamma}{1+c_{1}^{-1}}(1+c_{1}^{-1})<0.2\phi(1-\sigma)\gamma.
\end{equation}

From $\gamma<\frac{1-\sigma}{40000(1+c_{2}^{-1})L_{f}^{2}}$, and $\eta<\frac{(1-\sigma)^{2}\gamma}{40L_{f}}$, we have
\begin{equation}\label{omega2d}
\xi_{7}<40.00625(1+c_{2}^{-1})\gamma^{2}<\frac{(1-\sigma)^{2}}{320L_{f}^{2}}0.321(1-\sigma)\gamma=0.321\phi(1-\sigma)\gamma.
\end{equation}

From $\gamma<\frac{1-\sigma}{160(1+c_{1}^{-1})}$, and (\ref{omega2a})--(\ref{omega2d}), we know that (\ref{omega2}) holds. This completes the proof.

(iii) We third prove that the following inequality holds.
\begin{equation}\label{omega3}
0<\theta_{7}<1-0.44c_{1}(2c_{1}+1)<1.
\end{equation}

From $\gamma<\frac{c_{1}(1-\sigma)}{40C}$, we have
\begin{equation}\label{omega3a}
\xi_{8}<\frac{8}{1-\sigma}\frac{c_{1}(1-\sigma)C}{40C}<0.2c_{1}(2c_{1}+1).
\end{equation}

From $\gamma<\frac{c_{1}(1-\sigma)}{40C}$, we have
\begin{equation}\label{omega3b}
\phi\xi_{9}<\frac{(1-\sigma)^{2}}{320L_{f}^{2}}\frac{32L_{f}^{2}\gamma^{2}C}{(1-\sigma)\gamma}
<\frac{(1-\sigma)C}{10}\frac{c_{1}(1-\sigma)}{40C}<0.0025c_{1}(2c_{1}+1).
\end{equation}

From $\gamma<\frac{c_{1}}{8\sqrt{C}}$, we have
\begin{equation}\label{omega3c}
16\gamma^{2}(1+c_{1}^{-1})C<\frac{16c_{1}^{2}(1+c_{1}^{-1})C}{64C}<0.25c_{1}(2c_{1}+1).
\end{equation}

From $\gamma<\frac{c_{1}}{10L_{f}\sqrt{C(1+c_{2}^{-1})}}$, we have
\begin{align}\label{omega3d}
\xi_{11}<\frac{20(1+c_{2}^{-1})L_{f}^{2}c_{1}^{2}C}{100L_{f}^{2}C(1+c_{2}^{-1})}<0.1c_{1}(2c_{1}+1).
\end{align}

Considering that $\varphi_{X}\in(0,\frac{1}{r})$ and $\psi\in(0,1]$, we can ensure that $c_{1}=\frac{\varphi_{X}\psi r}{2}<\frac{1}{2}$. Consequently, we have $c_{1}(1+2c_{1})<1$.
Based on (\ref{omega3a})--(\ref{omega3d}), it becomes evident that (\ref{omega3}) is satisfied.
This completes the proof.

(iv) We fourth prove that the following inequality holds.
\begin{equation}\label{omega4}
0<\theta_{8}<1-0.77c_{2}(2c_{2}+1)<1.
\end{equation}


From $\gamma<\frac{c_{2}L_{f}^{2}}{C}$, we have
\begin{equation}\label{omega4a}
\phi\xi_{8}<\frac{(1-\sigma)^{2}}{320L_{f}^{2}}\frac{8c_{2}L_{f}^{2}C}{(1-\sigma)C}<0.025c_{2}(2c_{2}+1).
\end{equation}

From $\gamma<\frac{c_{2}}{10\sqrt{C}}$, we have
\begin{equation}\label{omega4b}
20\gamma^{2}(1+c_{2}^{-1})C<20\frac{c_{2}^{2}(1+c_{2}^{-1})C}{100C}<0.2c_{2}(1+2c_{2}).
\end{equation}

Considering that $\varphi_{Y}\in(0,\frac{1}{r})$ and $\psi\in(0,1]$, we can ensure that $c_{2}=\frac{\varphi_{Y}\psi r}{2}<\frac{1}{2}$. Consequently, we have $c_{2}(1+2c_{2})<1$.
From (\ref{omega4a})--(\ref{omega4b}), we know that (\ref{omega4}) holds. This completes the proof.

(v) From $0<\eta<\min\{\frac{1-\sigma}{160(1+c_{1}^{-1})},\frac{9}{40(4(1+c_{1}^{-1})+5(1+c_{2}^{-1}))}\}$, we have
\begin{equation}\label{omega5}
\theta_{2}>\eta(\frac{1}{4}-\frac{(1-\sigma)^{2}}{320L_{f}^{2}}\frac{8L_{f}^{2}}{(1-\sigma)}-\eta(4(1+c_{1}^{-1})+5(1+c_{2}^{-1}))>0.
\end{equation}

(vi) From $0<\eta<\frac{1}{2L_{f}}$, we have
\begin{equation}\label{omega6}
\theta_{9}=\frac{\eta}{4}(1-2\eta L_{f})>0.
\end{equation}

Based on (\ref{Lyapunovkp1}), (\ref{omega1}), (\ref{omega2}), (\ref{omega3}), (\ref{omega4}), (\ref{omega5}), (\ref{omega6}), we have
\begin{subequations}
\begin{align}
\bm{\mathrm{E}}_{\mathcal{C}}[U(k+1)]
\leq&\hspace{-0.3mm}\bm{\mathrm{E}}_{\mathcal{C}}[U(k)]-\theta_{3}\bm{\mathrm{E}}_{\mathcal{C}}[V(k)]-\theta_{2}\bm{\mathrm{E}}_{\mathcal{C}}[\|\bar{\bm{g}}^{0}(k)\|^{2}],\label{Lyapunovkp1b}\\
\leq&\hspace{-0.3mm}\bm{\mathrm{E}}_{\mathcal{C}}[U(k)]-\theta_{1}\bm{\mathrm{E}}_{\mathcal{C}}[\|\bm{X}(k)-\bar{\bm{X}}(k)\|^{2}+\|\bar{\bm{g}}^{0}(k)\|^{2}].\label{Lyapunovkp1c}
\end{align}
\end{subequations}


From (\ref{Lyapunovkp1c}), we have
\begin{equation}\label{Lyapunovkp1d}
\sum_{t=0}^{k}\bm{\mathrm{E}}_{\mathcal{C}}[\|\bm{X}(t)-\bar{\bm{X}}(t)\|^{2}+\|\bar{\bm{g}}^{0}(t)\|^{2}]\leq\frac{U(0)}{\theta_{1}},
\end{equation}
and
\begin{equation}\label{Lyapunovkp1e}
\bm{\mathrm{E}}_{\mathcal{C}}[n(F(\bar{X}(k))-F^{\star})]\leq\bm{\mathrm{E}}_{\mathcal{C}}[U(k)]<U(0).
\end{equation}

Based on (\ref{Lyapunovkp1d}) and (\ref{Lyapunovkp1e}), we know that (\ref{thm-1a}) and (\ref{thm-1b}) hold. This completes the proof.\\

\hspace{-3mm}\emph{C. Proof of Theorem 2}

From Assumption~\ref{assum-7}, we have
  \begin{equation}\label{barg0k}
  \|\bar{\bm{g}}^{0}(k)\|^{2}=n\|\nabla F(\bar{X}(k))\|^{2}\geq2\nu n(F(\bar{X}(k))-F^{\star}).
  \end{equation}

Based on (\ref{barg0k}), and (\ref{Lyapunovkp1b}), we have
\begin{align}\label{Lyapunovkp1f}
\bm{\mathrm{E}}_{\mathcal{C}}[U(k+1)]
&\leq\bm{\mathrm{E}}_{\mathcal{C}}[U(k)]-\theta_{3}\bm{\mathrm{E}}_{\mathcal{C}}[V(k)]-2\nu n\theta_{2}\bm{\mathrm{E}}_{\mathcal{C}}[n(F(\bar{X}(k))-F^{\star})]\notag\\
&\leq(1-\theta_{4})\bm{\mathrm{E}}_{\mathcal{C}}[U(k)]\leq(1-\theta_{4})^{k+1}U(0),
\end{align}
where $\theta_{4}=\min\{\theta_{3},2\nu\theta_{2}\}$.

From (\ref{omega2}), (\ref{omega3}), and (\ref{omega4}), we have $0<\theta_{4}\leq\theta_{3}<1$.
Based on (\ref{Lyapunovkp1f}), we know that (\ref{thm-2a}) holds.
This completes the proof.\\

\hspace{-3mm}\emph{D. Proof of Theorem~\ref{thm-3}}

We denote
  \begin{equation}\label{hatx}
  \hat{\bm{X}}(k)=\bm{A}(k)+\hat{\bm{Q}}^{X}(k),
  \end{equation}
  \begin{equation}\label{haty}
  \hat{\bm{Y}}(k)=\bm{C}(k)+\hat{\bm{Q}}^{Y}(k).
  \end{equation}

Then, the compact form of (\ref{Alg2g}) and (\ref{Alg2h}) respectively can be written as (\ref{reAlg1}), since $\bm{B}(k)=(\bm{\mathrm{I}}_{nd}-\bm{W})\bm{A}(k)$ and $\bm{D}(k)=(\bm{\mathrm{I}}_{nd}-\bm{W})\bm{C}(k)$.

From (\ref{Alg2i}), (\ref{Alg2k}), and (\ref{Assum1a}), we have
  \begin{align}\label{exkP1}
  &\bm{\mathrm{E}}_{\mathcal{C}}[\|\bm{E}^{X}(k+1)\|^{2}]\notag\\
  &\leq\bm{\mathrm{E}}_{\mathcal{C}}[C\|\varsigma\bm{E}^{X}(k)+\bm{X}(k)-\bm{A}(k)\|^{2}]\notag\\
  &\leq\bm{\mathrm{E}}_{\mathcal{C}}[2C\varsigma^{2}\|\bm{E}^{X}(k)\|^{2}+2C\|\bm{X}(k)-\bm{A}(k)\|^{2}],
  \end{align}
where the last inequality holds due to the Cauchy--Schwarz inequality.




Based on (\ref{reAlg1a}), (\ref{hatx}), and (\ref{exkP1}), we have
  \begin{align}\label{xhatx}
  &\hspace{-3mm}\bm{\mathrm{E}}_{\mathcal{C}}[\|\bm{X}(k)-\hat{\bm{X}}(k)\|^{2}]\notag\\
&\hspace{-3mm}=\bm{\mathrm{E}}_{\mathcal{C}}[\|\bm{E}^{X}(k+1)-\varsigma\bm{E}^{X}(k)\|^{2}]\notag\\
&\hspace{-3mm}\leq\bm{\mathrm{E}}_{\mathcal{C}}[2\|\bm{E}^{X}(k+1)\|^{2}+2\varsigma^{2}\|\bm{E}^{X}(k)\|^{2}]\notag\\
&\hspace{-3mm}\leq\bm{\mathrm{E}}_{\mathcal{C}}[2\varsigma^{2}(2C+1)\|\bm{E}^{X}(k)\|^{2}+4C\|\bm{X}(k)-\bm{A}(k)\|^{2}],
  \end{align}
where the first inequality holds due to the Cauchy--Schwarz inequality.

Then, from (\ref{ErrorxbarxnoE}) and (\ref{xhatx}), we have
\begin{align}\label{Errorxbarx1}
&\bm{\mathrm{E}}_{\mathcal{C}}[\|\bm{X}(k+1)-\bar{\bm{X}}(k+1)\|^{2}]\notag\\
&\leq\bm{\mathrm{E}}_{\mathcal{C}}[\frac{8\gamma}{1-\sigma}(2\varsigma^{2}(2C+1)\|\bm{E}^{X}(k)\|^{2}+4C\|\bm{X}(k)-\bm{A}(k)\|^{2})+\delta\|\bm{X}(k)-\bar{\bm{X}}(k)\|^{2}\notag\\
&\quad+\frac{2\eta^{2}}{\gamma(1-\sigma)}\|\bm{Y}(k)-\bar{\bm{Y}}(k)\|^{2}].
\end{align}


From (\ref{Alg2j}), (\ref{Alg2l}), and (\ref{Assum1a}), we have
  \begin{align}\label{eykP1}
  &\bm{\mathrm{E}}_{\mathcal{C}}[\|\bm{E}^{Y}(k+1)\|^{2}]\notag\\
  &\leq\bm{\mathrm{E}}_{\mathcal{C}}[C\|\varsigma\bm{E}^{Y}(k)+\bm{Y}(k)-\bm{C}(k)\|^{2}]\notag\\
  &\leq\bm{\mathrm{E}}_{\mathcal{C}}[2C\varsigma^{2}\|\bm{E}^{Y}(k)\|^{2}+2C\|\bm{Y}(k)-\bm{C}(k)\|^{2}],
  \end{align}
where the last inequality holds due to the Cauchy--Schwarz inequality.

Based on (\ref{reAlg1b}), (\ref{haty}), and (\ref{eykP1}), we have
  \begin{align}\label{yhaty}
  &\bm{\mathrm{E}}_{\mathcal{C}}[\|\bm{Y}(k)-\hat{\bm{Y}}(k)\|^{2}]\notag\\
  &=\bm{\mathrm{E}}_{\mathcal{C}}[\|\bm{E}^{Y}(k+1)-\varsigma\bm{E}^{Y}(k)\|^{2}]\notag\\
  &\leq\bm{\mathrm{E}}_{\mathcal{C}}[2\|\bm{E}^{Y}(k+1)\|^{2}+2\varsigma^{2}\|\bm{E}^{Y}(k)\|^{2}]\notag\\
  &\leq\bm{\mathrm{E}}_{\mathcal{C}}[2\varsigma^{2}(2C+1)\|\bm{E}^{Y}(k)\|^{2}+4C\|\bm{Y}(k)-\bm{C}(k)\|^{2}],
  \end{align}
where the first inequality holds due to the Cauchy--Schwarz inequality.

Then, from (\ref{Erroryg}), (\ref{xhatx}), and (\ref{yhaty}), we have
\begin{align}\label{Erroryy2}
&\bm{\mathrm{E}}_{\mathcal{C}}[\|\bm{Y}(k+1)-\bar{\bm{Y}}(k+1)\|^{2}]\notag\\
&\leq\bm{\mathrm{E}}_{\mathcal{C}}[\frac{8\gamma}{1-\sigma}(2\varsigma^{2}(2C+1)\|\bm{E}^{Y}(k)\|^{2}+4C\|\bm{Y}(k)-\bm{C}(k)\|^{2})+\delta\|\bm{Y}(k)-\bar{\bm{Y}}(k)\|^{2}\notag\\
&\quad+\frac{8L_{f}^{2}}{\gamma(1-\sigma)}(4\gamma^{2}(2\varsigma^{2}(2C+1)\|\bm{E}^{X}(k)\|^{2}+4C\|\bm{X}(k)-\bm{A}(k)\|^{2})\notag\\
&\quad+(8\gamma^{2}+2\eta^{2} L_{f}^{2})\|\bm{X}(k)-\bar{\bm{X}}(k)\|^{2}+\eta^{2}\|\bm{Y}(k)-\bar{\bm{Y}}(k)\|^{2}+\eta^{2}\|\bar{\bm{g}}^{0}(k)\|^{2})].
\end{align}



From (\ref{ExpErrorxa1}) and (\ref{xhatx}), we have
\begin{align}\label{ExpErrorxa2}
&\bm{\mathrm{E}}_{\mathcal{C}}[\|\bm{X}(k+1)-\bm{A}(k+1)\|^{2}]\notag\\
\leq&\bm{\mathrm{E}}_{\mathcal{C}}[4(1+c_{1}^{-1})(4\gamma^{2}\|\bm{X}(k)-\hat{\bm{X}}(k)\|^{2}+(8\gamma^{2}+2\eta^{2}L_{f}^{2})\|\bm{X}(k)-\bar{\bm{X}}(k)\|^{2}\notag\\
&+\eta^{2}\|\bm{Y}(k)-\bar{\bm{Y}}(k)\|^{2}+\eta^{2}\|\bar{\bm{g}}^{0}(k)\|^{2})+(1-c_{1}-2c_{1}^{2})\|\bm{X}(k)-\bm{A}(k)\|^{2}]\notag\\
\leq&\bm{\mathrm{E}}_{\mathcal{C}}[4(1+c_{1}^{-1})(4\gamma^{2}(2\varsigma^{2}(2C+1)\|\bm{E}^{X}(k)\|^{2}+4C\|\bm{X}(k)-\bm{A}(k)\|^{2})\notag\\
&+(8\gamma^{2}+2\eta^{2}L_{f}^{2})\|\bm{X}(k)-\bar{\bm{X}}(k)\|^{2}+\eta^{2}\|\bm{Y}(k)-\bar{\bm{Y}}(k)\|^{2}+\eta^{2}\|\bar{\bm{g}}^{0}(k)\|^{2})\notag\\
&+(1-c_{1}-2c_{1}^{2})\|\bm{X}(k)-\bm{A}(k)\|^{2}].
\end{align}



From (\ref{ExpErroryc2}), (\ref{xhatx}), and (\ref{yhaty}), we have
\begin{align}\label{ExpErroryc3}
&\bm{\mathrm{E}}_{\mathcal{C}}[\|\bm{Y}(k+1)-\bm{C}(k+1)\|^{2}]\notag\\
&\leq\bm{\mathrm{E}}_{\mathcal{C}}[5(1+c_{2}^{-1})(4\gamma^{2}(2\varsigma^{2}(2C+1)\|\bm{E}^{Y}(k)\|^{2}+4C\|\bm{Y}(k)-\bm{C}(k)\|^{2})\notag\\
&\quad+(8\gamma^{2}+2\eta^{2}L_{f}^{2})\|\bm{Y}(k)-\bar{\bm{Y}}(k)\|^{2}+L_{f}^{2}(4\gamma^{2}(2\varsigma^{2}(2C+1)\|\bm{E}^{X}(k)\|^{2}\notag\\
&\quad+4C\|\bm{X}(k)-\bm{A}(k)\|^{2})+(8\gamma^{2}+2\eta^{2}L_{f}^{2})\|\bm{X}(k)-\bar{\bm{X}}(k)\|^{2}+\eta^{2}\|\bar{\bm{g}}^{0}(k)\|^{2}))\notag\\
&\quad+(1-c_{2}-2c_{2}^{2})\|\bm{Y}(k)-\bm{C}(k)\|^{2}].
\end{align}

From (\ref{exkP1}), (\ref{Errorxbarx1}), (\ref{eykP1}), (\ref{Erroryy2}), (\ref{ExpErrorxa2}), and (\ref{ExpErroryc3}), we have
\begin{align}\label{Lyapunovkp2}
&\bm{\mathrm{E}}_{\mathcal{C}}[\hat{U}(k+1)]\notag\\
&\leq\bm{\mathrm{E}}_{\mathcal{C}}[\theta_{5}\|\bm{X}(k)-\bar{\bm{X}}(k)\|^{2}+\theta_{6}\|\bm{Y}(k)-\bar{\bm{Y}}(k)\|^{2}+\hat{\theta}_{4}\|\bm{X}(k)-\bm{A}(k)\|^{2}+\hat{\theta}_{5}\|\bm{Y}(k)-\bm{C}(k)\|^{2}\notag\\
&\quad+\hat{\theta}_{6}\|\bm{E}^{X}(k)\|^{2}+\hat{\theta}_{7}\|\bm{E}^{Y}(k)\|^{2}+n(F(\bar{X}(k))-F^{\star})-\theta_{2}\|\bar{\bm{g}}^{0}(k)\|^{2}-\theta_{9}\|\bar{\bm{g}}(k)\|^{2}],
\end{align}
where
\begin{align*}
\hat{\theta}_{4}&=4\xi_{8}+4\phi\xi_{9}+\hat{\xi}_{3}+4\xi_{11}+2C\hat{\phi},~\hat{\theta}_{5}=4\phi\xi_{8}+\hat{\xi}_{5}+2C\hat{\phi},\notag\\
\hat{\theta}_{6}&=\hat{\xi}_{1}+\phi\hat{\xi}_{2}+\hat{\xi}_{4}+\hat{\xi}_{6}+2C\hat{\phi}\varsigma^{2},~\hat{\theta}_{7}=\phi\hat{\xi}_{1}+\hat{\xi}_{7}+2C\hat{\phi}\varsigma^{2},\notag\\
\hat{\xi}_{1}&=\frac{8\gamma}{1-\sigma}2\varsigma^{2}(2C+1),~\hat{\xi}_{2}=\frac{8L_{f}^{2}}{(1-\sigma)\gamma}4\gamma^{2}(2\varsigma^{2}(2C+1)),\\
\hat{\xi}_{3}&=64\gamma^{2}(1+c_{1}^{-1})C+(1-c_{1}-2c_{1}^{2}),~\hat{\xi}_{4}=16\gamma^{2}(1+c_{1}^{-1})(2\varsigma^{2}(2C+1)),\\
\hat{\xi}_{5}&=80\gamma^{2}(1+c_{2}^{-1})C+(1-c_{2}-2c_{2}^{2}),~\hat{\xi}_{6}=L_{f}^{2}\hat{\xi}_{7},~\hat{\xi}_{7}=5(1+c_{2}^{-1})4\gamma^{2}(2\varsigma^{2}(2C+1)).
\end{align*}

Next, we prove that $\hat\theta_{4},~\hat\theta_{5},~\hat\theta_{6},~\hat\theta_{7}\in(0,1)$.

(i) We third prove that the following inequality holds.
\begin{equation}\label{hatomega1}
0<\hat{\theta}_{4}<1-0.24c_{1}(2c_{1}+1)<1.
\end{equation}

From $\gamma<\frac{c_{1}(1-\sigma)}{160C}$, we have
\begin{equation}\label{hatomega1a}
4\xi_{8}<\frac{32}{1-\sigma}\frac{c_{1}(1-\sigma)C}{160C}<0.2c_{1}(2c_{1}+1).
\end{equation}

From $\gamma<\frac{c_{1}(1-\sigma)}{160C}$, we have
\begin{equation}\label{hatomega1b}
4\phi\xi_{9}<4\frac{(1-\sigma)^{2}}{320L_{f}^{2}}\frac{32L_{f}^{2}\gamma^{2}C}{(1-\sigma)\gamma}
<\frac{4(1-\sigma)C}{10}\frac{c_{1}(1-\sigma)}{160C}<0.0025c_{1}(2c_{1}+1).
\end{equation}

From $\gamma<\frac{c_{1}}{16\sqrt{C}}$, we have
\begin{equation}\label{hatomega1c}
64\gamma^{2}(1+c_{1}^{-1})C<\frac{64c_{1}^{2}(1+c_{1}^{-1})C}{256C}<0.25c_{1}(2c_{1}+1).
\end{equation}

From $\gamma<\frac{c_{1}}{20L_{f}\sqrt{C(1+c_{2}^{-1})}}$, we have
\begin{equation}\label{hatomega1d}
4\xi_{11}<\frac{20(1+c_{2}^{-1})L_{f}^{2}c_{1}^{2}C}{100L_{f}^{2}C(1+c_{2}^{-1})}<0.1c_{2}(2c_{2}+1).
\end{equation}

From $\hat{\phi}\leq\frac{0.1}{C}c_{1}(2c_{1}+1)$, we have
\begin{equation}\label{hatomega1e}
2C\hat{\phi}\leq0.2c_{1}(2c_{1}+1).
\end{equation}

Considering that $\varphi_{X}\in(0,\frac{1}{r})$ and $\psi\in(0,1]$, we can ensure that $c_{1}=\frac{\varphi_{X}\psi r}{2}<\frac{1}{2}$. Consequently, we have $c_{1}(1+2c_{1})<1$.
From (\ref{hatomega1a})--(\ref{hatomega1e}), we know that (\ref{hatomega1}) holds. This completes the proof.

(ii) We fourth prove that the following inequality holds.
\begin{equation}\label{hatomega2}
0<\hat{\theta}_{5}<1-0.57c_{2}(2c_{2}+1)<1.
\end{equation}

From $\gamma<\frac{c_{2}L_{f}^{2}}{4C}$, we have
\begin{equation}\label{hatomega2a}
4\phi\xi_{8}<\frac{4(1-\sigma)^{2}}{320L_{f}^{2}}\frac{8c_{2}L_{f}^{2}C}{4(1-\sigma)C}<0.025c_{2}(2c_{2}+1).
\end{equation}

From $\gamma<\frac{c_{2}}{20\sqrt{C}}$, we have
\begin{equation}\label{hatomega2b}
80\gamma^{2}(1+c_{2}^{-1})C<80\frac{c_{2}^{2}(1+c_{2}^{-1})C}{400C}<0.2c_{2}(1+2c_{2}).
\end{equation}

From $\hat{\phi}\leq\frac{0.1}{C}c_{2}(2c_{2}+1)$, we have
\begin{equation}\label{hatomega2c}
2C\hat{\phi}\leq0.2c_{2}(2c_{2}+1).
\end{equation}

Considering that $\varphi_{Y}\in(0,\frac{1}{r})$ and $\psi\in(0,1]$, we can ensure that $c_{2}=\frac{\varphi_{Y}\psi r}{2}<\frac{1}{2}$. Consequently, we have $c_{2}(1+2c_{2})<1$.
From (\ref{hatomega2a})--(\ref{hatomega2c}), we know that (\ref{hatomega2}) holds. This completes the proof.

(iii) We third prove that the following inequality holds.
\begin{equation}\label{hatomega3}
0<\hat{\theta}_{6}<0.75\hat{\phi}.
\end{equation}

From $\varsigma<\frac{1}{2\sqrt{C}}$, we have
\begin{equation}\label{hatomega3a}
2C\hat{\phi}\varsigma^{2}<0.5\hat{\phi}.
\end{equation}

From $\varsigma<\frac{1}{\sqrt{2C+1}}$, $\gamma<\frac{1}{4(1+c_{1}^{-1})+5(1+c_{2}^{-1})L_{f}^{2}}$, and $\gamma<\frac{(1-\sigma)\hat{\phi}}{4(16(1+4\phi L_{f}^{2})+8(1-\sigma))}$,
we have
\begin{align}\label{hatomega3b}
&\hat{\xi}_{1}+\phi\hat{\xi}_{2}+\hat{\xi}_{4}+\hat{\xi}_{6}\notag\\
&=\gamma(\frac{16(1+4\phi L_{f}^{2})}{1-\sigma}+8\gamma(4(1+c_{1}^{-1})+5(1+c_{2}^{-1})L_{f}^{2}))\varsigma^{2}(2C+1)\notag\\
&<\gamma(\frac{16(1+4\phi L_{f}^{2})}{1-\sigma}+8)<0.25\hat{\phi}.
\end{align}

From (\ref{hatomega3a}), and (\ref{hatomega3b}), we know that (\ref{hatomega3}) holds.

(iv) We fourth prove that the following inequality holds.
\begin{equation}\label{hatomega4}
0<\hat{\theta}_{7}<0.75\hat{\phi}.
\end{equation}

From $\varsigma<\frac{1}{\sqrt{2C+1}}$, $\gamma<\frac{1}{5(1+c_{2}^{-1})}$, and $\gamma<\frac{(1-\sigma)\hat{\phi}}{32(2\phi+(1-\sigma))}$,
we have
\begin{equation}\label{hatomega4a}
\phi\hat{\xi}_{1}+\hat{\xi}_{7}=\gamma(\frac{16\phi}{1-\sigma}+40\gamma(1+c_{2}^{-1}))\varsigma^{2}(2C+1)<\gamma(\frac{16\phi}{1-\sigma}+8)<0.25\hat{\phi}.
\end{equation}

From (\ref{hatomega3a}), and (\ref{hatomega4a}), we know that (\ref{hatomega4}) holds.

%
%

Based on (\ref{omega1}), (\ref{omega2}), (\ref{omega5}), (\ref{omega6}), (\ref{Lyapunovkp2}), (\ref{hatomega1}), (\ref{hatomega2}), (\ref{hatomega3}), and (\ref{hatomega4}), we have
\begin{subequations}
\begin{align}
&\bm{\mathrm{E}}_{\mathcal{C}}[\hat{U}(k+1)]\notag\\
\leq&\hspace{-0.3mm}\bm{\mathrm{E}}_{\mathcal{C}}[\hat{U}(k)]-\hat{\theta}_{2}\bm{\mathrm{E}}_{\mathcal{C}}[\hat{V}(k)]-\theta_{2}\bm{\mathrm{E}}_{\mathcal{C}}[\|\bar{\bm{g}}^{0}(k)\|^{2}],\label{Lyapunovkp2b}\\
\leq&\hspace{-0.3mm}\bm{\mathrm{E}}_{\mathcal{C}}[\hat{U}(k)]-\hat{\theta}_{1}\bm{\mathrm{E}}_{\mathcal{C}}[\|\bm{X}(k)-\bar{\bm{X}}(k)\|^{2}+\|\bar{\bm{g}}^{0}(k)\|^{2}].\label{Lyapunovkp2c}
\end{align}
\end{subequations}


From (\ref{Lyapunovkp2c}), we have
\begin{equation}\label{Lyapunovkp2d}
\sum_{t=0}^{k}\bm{\mathrm{E}}_{\mathcal{C}}[\|\bm{X}(t)-\bar{\bm{X}}(t)\|^{2}+\|\bar{\bm{g}}^{0}(t)\|^{2}]\leq\frac{\hat{U}(0)}{\hat{\theta}_{1}},
\end{equation}
and
\begin{equation}\label{Lyapunovkp2e}
\bm{\mathrm{E}}_{\mathcal{C}}[n(F(\bar{X}(k))-F^{\star})]\leq\bm{\mathrm{E}}_{\mathcal{C}}[\hat{U}(k)]<\hat{U}(0).
\end{equation}

Based on (\ref{Lyapunovkp2d}) and (\ref{Lyapunovkp2e}), we know that (\ref{thm-3a}) and (\ref{thm-3b}) hold. This completes the proof.\\

\hspace{-3mm}\emph{E. Proof of Theorem 4}


Based on (\ref{barg0k}), and (\ref{Lyapunovkp2b}), we have
\begin{align}\label{Lyapunovkp2f}
\bm{\mathrm{E}}_{\mathcal{C}}[U(k+1)]
&<\bm{\mathrm{E}}_{\mathcal{C}}[U(k)]-\hat{\theta}_{2}\bm{\mathrm{E}}_{\mathcal{C}}[V(k)]-2\nu n\theta_{2}\bm{\mathrm{E}}_{\mathcal{C}}[n(F(\bar{X}(k))-F^{\star})]\notag\\
&\leq(1-\hat{\theta}_{3})\bm{\mathrm{E}}_{\mathcal{C}}[U(k)]\leq(1-\hat{\theta}_{3})^{k+1}U(0),
\end{align}
where $\hat{\theta}_{3}=\min\{\hat{\theta}_{2},2\nu\theta_{2}\}$.

Based on (\ref{omega2}), (\ref{hatomega1}), and (\ref{hatomega2}), we have $0<\hat\theta_{3}\leq\hat\theta_{2}<1$.
Based on (\ref{Lyapunovkp2f}), we know that (\ref{thm-4a}) holds.
This completes the proof.\\

\hspace{-3mm}\emph{F. Proof of Theorem~\ref{thm-5}}

From Algorithm~\ref{alg:3}, it is straightforward to check that
\begin{align*}
V_{i}(k)&=\hat{X}_{i}(k)-\sum_{j=1}^{n}W_{ij}\hat{X}_{j}(k), \notag\\
Z_{i}(k)&=\hat{Y}_{i}(k)-\sum_{j=1}^{n}W_{ij}\hat{Y}_{j}(k).
\end{align*}

Then, noting that (\ref{Alg3e}) and (\ref{Alg3f}) can respectively be rewritten the compact form as (\ref{reAlg1a}) and (\ref{reAlg1b}), we know that (\ref{ErrorxbarxnoE}) and (\ref{ErrorygnoE}) still hold.

%

From \cite[Equation(5.4.21)]{Horn2012Matrix}, we have
\begin{align}\label{xminushatx}
\|\bm{X}(k)-\hat{\bm{X}}(k)\|^2&=\sum_{i=1}^n\|X_{i}(k)-\hat{X}_{i}(k)\|^2\notag\\
&\leq\sum_{i=1}^n\tilde{d}^2\|X_{i}(k)-\hat{X}_{i}(k)\|_p^2\notag\\
&\leq n\tilde{d}^2\max_{i\in\mathcal{V}}\|X_{i}(k)-\hat{X}_{i}(k)\|_p^2.
\end{align}

Then, from (\ref{Alg3a}), (\ref{Alg3g}), and (\ref{Assum2}), we have
\begin{align}\label{xhatxp}
&\bm{\mathrm{E}}_{\mathcal{C}}[\|X_{i}(k)-\hat{X}_{i}(k)\|_p^2]\notag\\
&=\bm{\mathrm{E}}_{\mathcal{C}}[\|X_{i}(k)-\hat{X}_{i}(k-1)-s(k)\mathcal{C}((X_{i}(k)-\hat{X}_{i}(k-1))/s(k))\|_p^2]\notag\\
&=\bm{\mathrm{E}}_{\mathcal{C}}[s(k)^{2}\|(X_{i}(k)-\hat{X}_{i}(k-1))/s(k)-\mathcal{C}((X_{i}(k)-\hat{X}_{i}(k-1))/s(k))\|_p^2]\notag\\
&\leq Cs^{2}(k)=C s^{2}(0)\mu^{2k}.
\end{align}

From (\ref{ErrorxbarxnoE}), (\ref{xminushatx}), and (\ref{xhatxp}), we have
\begin{align}\label{ErrorxbarxA3b}
&\bm{\mathrm{E}}_{\mathcal{C}}[\|\bm{X}(k+1)-\bar{\bm{X}}(k+1)\|^{2}]\notag\\
&\leq\bm{\mathrm{E}}_{\mathcal{C}}[\delta\|\bm{X}(k)-\bar{\bm{X}}(k)\|^{2}+\frac{2\eta^{2}}{\gamma(1-\sigma)}\|\bm{Y}(k)-\bar{\bm{Y}}(k)\|^{2}+n\tilde{d}^2\frac{8\gamma}{1-\sigma}Cs^{2}(0)\mu^{2k}].
\end{align}


From \cite[Equation(5.4.21)]{Horn2012Matrix}, we have
\begin{equation}\label{xyminushaty}
\|\bm{Y}(k)-\hat{\bm{Y}}(k)\|^2\leq n\tilde{d}^2\max_{i\in\mathcal{V}}\|Y_{i}(k)-\hat{Y}_{i}(k)\|_p^2.
\end{equation}

Then, from (\ref{Alg3b}), (\ref{Alg3h}), and (\ref{Assum2}), we have
\begin{equation}\label{yhatyp}
\bm{\mathrm{E}}_{\mathcal{C}}[\|Y_{i}(k)-\hat{Y}_{i}(k)\|_p^2]\leq Cs^{2}(0)\mu^{2k}.
\end{equation}

From (\ref{Erroryg}), (\ref{xyminushaty}), and (\ref{yhatyp}), we have
\begin{align}\label{ErrorybaryA3b}
&\bm{\mathrm{E}}_{\mathcal{C}}[\|\bm{Y}(k+1)-\bar{\bm{Y}}(k+1)\|^{2}]\notag\\
&\leq\bm{\mathrm{E}}_{\mathcal{C}}[(\delta+\frac{8L_{f}^{2}}{\gamma(1-\sigma)}\eta^{2})\|\bm{Y}(k)-\bar{\bm{Y}}(k)\|^{2}+\frac{8L_{f}^{2}}{\gamma(1-\sigma)}(4\gamma^{2}Cn\tilde{d}^2s^{2}(0)\mu^{2k}\notag\\
&\quad+(8\gamma^{2}+2\eta^{2} L_{f}^{2})\|\bm{X}(k)-\bar{\bm{X}}(k)\|^{2}+\eta^{2}\|\bar{\bm{g}}^{0}(k)\|^{2})+\frac{8\gamma}{1-\sigma}Cn\tilde{d}^2s^{2}(0)\mu^{2k}].
\end{align}


From (\ref{ErrorxbarxnoE}), (\ref{ErrorygnoE}), (\ref{ErrorxbarxA3b}), (\ref{ErrorybaryA3b}), we have
\begin{align}\label{Lyapunovkp3}
&\bm{\mathrm{E}}_{\mathcal{C}}[\breve{U}(k+1)]\notag\\
&\leq\bm{\mathrm{E}}_{\mathcal{C}}[\breve\theta_{12}\|\bm{X}(k)-\bar{\bm{X}}(k)\|^{2}+\breve\theta_{13}\|\bm{Y}(k)-\bar{\bm{Y}}(k)\|^{2}+n(F(\bar{X}(k))-F^{\star})\notag\\
&\quad-\breve\theta_{4}\|\bar{\bm{g}}^{0}(k)\|^{2}-\theta_{9}\|\bar{\bm{g}}(k)\|^{2}+2n\tilde{d}^2\xi_{8}s^{2}(k)(1+2L_{f}^{2})],
\end{align}
where
$$
\breve\theta_{12}=\delta+\phi\xi_{1}+\frac{\eta L_{f}^{2}}{2},~\breve\theta_{13}=\xi_{4}+\phi\xi_{5}.
$$

From (\ref{omega1a}), (\ref{omega1d}), (\ref{omega2a}), (\ref{omega2b}), and (\ref{omega5}), we have $\breve\theta_{4}>0$, and
\begin{equation}\label{breveomeg1}
0<\breve\theta_{12}<1-0.59(1-\sigma)\gamma<1,
\end{equation}
\begin{equation}\label{breveomeg2}
0<\breve\theta_{13}<1-0.59\phi(1-\sigma)\gamma<1.
\end{equation}

Based on (\ref{omega5}), (\ref{omega6}), (\ref{Lyapunovkp3}), (\ref{breveomeg1}), and (\ref{breveomeg2}), we have
\begin{subequations}
\begin{align}
&\bm{\mathrm{E}}_{\mathcal{C}}[\breve{U}(k+1)]\notag\\
\leq&\bm{\mathrm{E}}_{\mathcal{C}}[\breve{U}(k)]-\breve\theta_{3}\bm{\mathrm{E}}_{\mathcal{C}}[\breve{V}(k)]-\breve\theta_{4}\bm{\mathrm{E}}_{\mathcal{C}}[\|\bar{\bm{g}}^{0}(k)\|^{2}+2n\tilde{d}^2\xi_{8}s^{2}(k)(1+2L_{f}^{2})\notag\\
\leq&\bm{\mathrm{E}}_{\mathcal{C}}[\breve{U}(k)]-\breve\theta_{1}\bm{\mathrm{E}}_{\mathcal{C}}[\|\bm{X}(k)-\bar{\bm{X}}(k)\|^{2}+\|\bar{\bm{g}}^{0}(k)\|^{2}]+2n\tilde{d}^2\xi_{8}s^{2}(k)(1+2L_{f}^{2})\notag\\
\leq&\breve{U}(0)-\sum_{t=0}^{k}\breve\theta_{1}\bm{\mathrm{E}}_{\mathcal{C}}[\|\bm{X}(t)-\bar{\bm{X}}(t)\|^{2}+\|\bar{\bm{g}}^{0}(t)\|^{2}]+\sum_{t=0}^{k}2n\tilde{d}^2\xi_{8}s^{2}(t)(1+2L_{f}^{2}),\label{Lyapunovkp3b}\\
\leq&\breve{U}(0)-\sum_{t=0}^{k}\breve\theta_{1}\bm{\mathrm{E}}_{\mathcal{C}}[\|\bm{X}(t)-\bar{\bm{X}}(t)\|^{2}+\|\bar{\bm{g}}^{0}(t)\|^{2}]+\frac{\breve\theta_{2}}{1-\mu^{2}}.\label{Lyapunovkp3c}
\end{align}
\end{subequations}


From (\ref{Lyapunovkp3c}), we have
\begin{equation}\label{Lyapunovkp3d}
\sum_{t=0}^{k}\bm{\mathrm{E}}_{\mathcal{C}}[\|\bm{X}(t)-\bar{\bm{X}}(t)\|^{2}+\|\bar{\bm{g}}^{0}(t)\|^{2}]\leq\frac{\breve{U}(0)+\frac{\breve\theta_{2}}{1-\mu^{2}}}{\breve\theta_{1}},
\end{equation}
and
\begin{equation}\label{Lyapunovkp3e}
\bm{\mathrm{E}}_{\mathcal{C}}[n(F(\bar{X}(k))-F^{\star})]\leq\bm{\mathrm{E}}_{\mathcal{C}}[\breve{U}(k)]<\breve{U}(0)+\frac{\breve\theta_{2}}{1-\mu^{2}}.
\end{equation}

Based on (\ref{Lyapunovkp3d}) and (\ref{Lyapunovkp3e}), we know that (\ref{thm-5a}) and (\ref{thm-5b}) hold. This completes the proof.\\

\hspace{-3mm}\emph{G. Proof of Theorem 6}

Based on (\ref{barg0k}), and (\ref{Lyapunovkp3b}), we have
\begin{align}\label{Lyapunovkp3f}
&\bm{\mathrm{E}}_{\mathcal{C}}[\breve{U}(k+1)]\notag\\
&<\bm{\mathrm{E}}_{\mathcal{C}}[\breve{U}(k)]-\breve\theta_{1}\bm{\mathrm{E}}_{\mathcal{C}}[\breve{V}(k)]-2\nu n\theta_{2}\bm{\mathrm{E}}_{\mathcal{C}}[n(F(\bar{X}(k))-F^{\star})]+2n\tilde{d}^2\xi_{8}s^{2}(k)(1+2L_{f}^{2})\notag\\
&\leq(1-\breve\theta_{7})\bm{\mathrm{E}}_{\mathcal{C}}[\breve{U}(k)]+2n\tilde{d}^2\xi_{8}s^{2}(k)(1+2L_{f}^{2})\notag\\
&\leq(1-\breve\theta_{7})^{k+1}\breve{U}(0)+\sum_{t=0}^{k}(1-\breve\theta_{7})^{t}\breve\theta_{8}s^{2}(k-t).
\end{align}

If $1-\breve\theta_{7}<\mu^{2}$, then
\begin{align}\label{Lyapunovkp3fa}
\bm{\mathrm{E}}_{\mathcal{C}}[\breve{U}(k+1)]
&\leq(1-\breve\theta_{7})^{k+1}\breve{U}(0)+\sum_{t=0}^{k}\breve\theta_{9}^{t}\breve\theta_{8}s^{2}(0)\mu^{2k}\notag\\
&=(1-\breve\theta_{7})^{k+1}\breve{U}(0)+\frac{1-\breve\theta_{9}^{k}}{1-\breve\theta_{9}}\breve\theta_{8}s^{2}(0)\mu^{2k}\notag\\
&\leq\mu^{2(k+1)}(\breve{U}(0)+\frac{\breve\theta_{8}s^{2}(0)}{(1-\breve\theta_{9})\mu^{2}}).
\end{align}

If $1-\breve\theta_{7}>\mu^{2}$, then
\begin{align}\label{Lyapunovkp3fb}
\bm{\mathrm{E}}_{\mathcal{C}}[\breve{U}(k+1)]
&\leq(1-\breve\theta_{7})^{k+1}\breve{U}(0)+\sum_{t=0}^{k}\breve\theta_{10}^{t}\breve\theta_{8}s^{2}(0)(1-\breve\theta_{7})^{k}\notag\\
&=(1-\breve\theta_{7})^{k+1}\breve{U}(0)+\frac{1-\breve\theta_{10}^{k}}{1-\breve\theta_{10}}\breve\theta_{8}s^{2}(0)(1-\breve\theta_{7})^{k}\notag\\
&\leq(1-\breve\theta_{7})^{k+1}(\breve{U}(0)+\frac{\breve\theta_{8}s^{2}(0)}{(1-\breve\theta_{10})(1-\breve\theta_{7})}).
\end{align}

If $1-\breve\theta_{7}=\mu^{2}$, then for any $\varpi\in(\mu^{2},1)$, we have
\begin{align}\label{Lyapunovkp3fc}
\bm{\mathrm{E}}_{\mathcal{C}}[\breve{U}(k+1)]
&\leq\varpi^{k+1}\breve{U}(0)+\sum_{t=0}^{k}\breve\theta_{11}^{t}\breve\theta_{8}s^{2}(0)\varpi^{k}\notag\\
&=\varpi^{k+1}\breve{U}(0)+\frac{1-\breve\theta_{11}^{k}}{1-\breve\theta_{11}}\breve\theta_{8}s^{2}(0)\varpi^{k}\notag\\
&\leq\varpi^{k+1}(\breve{U}(0)+\frac{\breve\theta_{8}s^{2}(0)}{(1-\breve\theta_{11})\varpi}).
\end{align}

From $1>\breve\theta_{3}>0$ and (\ref{omega5}), we know that $\breve\theta_{1}>0$,
and
\begin{align}\label{breve_omega7}
0<\breve\theta_{7}\leq\breve\theta_{1}\leq\breve\theta_{3}<1.
\end{align}

From (\ref{breve_omega7}) and (\ref{omega5}), we have
\begin{align}\label{breve_omega8}
0<\breve\theta_{5}\leq1-\mu^{2}<1.
\end{align}

Then, based on (\ref{Lyapunovkp3fa})--(\ref{breve_omega8}), we know that (\ref{thm-6a}) holds. This completes the proof.\\

\hspace{-3mm}\emph{H. Proof of Theorem 7}

We use mathematical induction to prove the following inequalities hold.
  \begin{subequations}\label{Induction}
  \begin{align}
  &\tilde U(k)\leq\tilde{\xi}_{5}s^{2}(k),\label{Inductiona}\\
  &\max_{i\in\mathcal{V}}\|X_{i}(k)-\hat{X}_{i}(k-1)\|_{p}^{2}\leq s^{2}(k),\label{Inductionb}\\
  &\max_{i\in\mathcal{V}}\|Y_{i}(k)-\hat{Y}_{i}(k-1)\|_{p}^{2}\leq s^{2}(k).\label{Inductionc}
  \end{align}
  \end{subequations}

For $k = 0$, from $s(0)\geq\sqrt{\frac{\tilde U(0)}{\tilde{\xi}_{5}}}$, we have
  \begin{equation}\label{Induction0a}
  \tilde U(0)\leq\tilde{\xi}_{5}s^{2}(0),
  \end{equation}
where the first inequality holds due to (\ref{PL}).

From $s(0)\geq\max\{\max_{i\in\mathcal{V}}\|X_{i}(0)\|,\max_{i\in\mathcal{V}}\|Y_{i}(0)\|\}$, $\hat{X}_{i}(-1)=\bm{0}_{d}$, $\hat{Y}_{i}(-1)=\bm{0}_{d}$, and \cite[Equation(5.4.21)]{Horn2012Matrix}, we have
  \begin{subequations}\label{Induction0}
  \begin{align}
  \max_{i\in\mathcal{V}}\|X_{i}(0)-\hat{X}_{i}(-1)\|_{p}^{2}\leq s^{2}(0),\label{Induction0b}\\
  \max_{i\in\mathcal{V}}\|Y_{i}(0)-\hat{Y}_{i}(-1)\|_{p}^{2}\leq s^{2}(0).\label{Induction0c}
  \end{align}
  \end{subequations}

Therefore, from (\ref{Induction0a}) and (\ref{Induction0}), we know that (\ref{Induction}) holds at $k=0$.

Suppose that (\ref{Induction}) holds at $k$.
We next show that (\ref{Induction}) holds at $k+1$.

For $k = k+1$, we have
  \begin{align}\label{xhatxD}
  &\|X_{i}(k+1)-\hat{X}_{i}(k)\|_{p}^{2}\notag\\
  &=\|X_{i}(k+1)-X_{i}(k)+X_{i}(k)-\hat{X}_{i}(k)\|_{p}^{2}\notag\\
  &\leq(\|X_{i}(k+1)-X_{i}(k)\|_{p}+\|X_{i}(k)-\hat{X}_{i}(k)\|_{p})^{2}\notag\\
  &\leq(1+\varphi^{-1})\|X_{i}(k+1)-X_{i}(k)\|_{p}^{2}+(1+\varphi)\|X_{i}(k)-\hat{X}_{i}(k)\|_{p}^{2}\notag\\
  &\leq(1+\varphi^{-1})\hat{d}^{2}\|X_{i}(k+1)-X_{i}(k)\|^{2}+(1+\varphi)\|X_{i}(k)-\hat{X}_{i}(k)\|_{p}^{2}\notag\\
  &\leq(1+\varphi^{-1})\hat{d}^{2}\|\bm{X}(k+1)-\bm{X}(k)\|^{2}+(1+\varphi)\|X_{i}(k)-\hat{X}_{i}(k)\|_{p}^{2},
  \end{align}
where the first inequality holds due to the Minkowski inequality; the second inequality holds due to the Cauchy--Schwarz inequality; the third inequality holds due to
\cite[Equation(5.4.21)]{Horn2012Matrix}.

Similarly, we have
  \begin{align}\label{yhatyD}
  &\|Y_{i}(k+1)-\hat{Y}_{i}(k)\|_{p}^{2}\notag\\
  &\leq(1+\varphi^{-1})\hat{d}^{2}\|\bm{Y}(k+1)-\bm{Y}(k)\|^{2}+(1+\varphi)\|Y_{i}(k)-\hat{Y}_{i}(k)\|_{p}^{2}.
  \end{align}

Based on (\ref{Errorxx}) and (\ref{xhatxD}), we have
  \begin{align}\label{xhatxD2}
  &\|X_{i}(k+1)-\hat{X}_{i}(k)\|_{p}^{2}\notag\\
  &\leq4(1+\varphi^{-1})\hat{d}^{2}(4\gamma^{2}\|\bm{X}(k)-\bm{\hat{X}}(k)\|^{2}+(8\gamma^{2}+2\eta^{2}L_{f}^{2})\|\bm{X}(k)-\bar{\bm{X}}(k)\|^{2}\notag\\
  &\quad+\eta^{2}\|\bm{Y}(k)-\bar{\bm{Y}}(k)\|^{2}+\eta^{2}\|\bar{\bm{g}}^{0}(k)\|^{2})+(1+\varphi)\|X_{i}(k)-\hat{X}_{i}(k)\|_{p}^{2}.
  \end{align}

Similarly, based on (\ref{Erroryy}) and (\ref{yhatyD}), we have
  \begin{align}\label{yhatyD2}
  &\|Y_{i}(k+1)-\hat{Y}_{i}(k)\|_{p}^{2}\notag\\
  &\leq5(1+\varphi^{-1})\hat{d}^{2}(4\gamma^{2}\|\bm{Y}(k)-\bm{\hat{Y}}(k)\|^{2}+(8\gamma^{2}+2\eta^{2}L_{f}^{2})\|\bm{Y}(k)-\bar{\bm{Y}}(k)\|^{2}\notag\\
  &\quad+L_{f}^{2}(4\gamma^{2}\|\bm{X}(k)-\bm{\hat{X}}(k)\|^{2}+(8\gamma^{2}+2\eta^{2}L_{f}^{2})\|\bm{X}(k)-\bar{\bm{X}}(k)\|^{2}+\eta^{2}\|\bar{\bm{g}}^{0}(k)\|^{2}))\notag\\
  &\quad+(1+\varphi)\|Y_{i}(k)-\hat{Y}_{i}(k)\|_{p}^{2}.
  \end{align}

From (\ref{OptimalF2}), we know that
\begin{equation}\label{OptimalF2D}
\|\bar{\bm{g}}^{0}(k)\|^{2}\leq\frac{4}{\eta}(1+\frac{\eta L_{f}^{2}}{2})(n(F(\bar{X}(k))-F^{\star})+\|\bm{X}(k)-\bar{\bm{X}}(k)\|^{2}).
\end{equation}

Then, based on (\ref{xminushatx}), (\ref{xhatxD2}), and (\ref{OptimalF2D}), we have
  \begin{equation}\label{xhatxD3}
  \max_{i\in\mathcal{V}}\|X_{i}(k+1)-\hat{X}_{i}(k)\|_{p}^{2}\leq(1+\varphi+\tilde{\xi}_{8}\gamma^{2})\|X_{i}(k)-\hat{X}_{i}(k)\|_{p}^{2}+\tilde{\xi}_{10}(1+\varphi^{-1})\tilde U(k),
  \end{equation}
where $\tilde{\xi}_{10}=4\hat{d}^{2}(8\gamma^{2}+\eta^{2}(1+4L_{f}^{2})+4\eta)$.

From (\ref{xyminushaty}), (\ref{yhatyD2}), and (\ref{OptimalF2D}), we have
  \begin{align}\label{yhatyD3}
  &\max_{i\in\mathcal{V}}\|Y_{i}(k+1)-\hat{Y}_{i}(k)\|_{p}^{2}\notag\\
  &\leq(1+\varphi+\tilde{\xi}_{9}\gamma^{2})\|Y_{i}(k)-\hat{Y}_{i}(k)\|_{p}^{2}+\tilde{\xi}_{9}L_{f}^{2}\eta\|X_{i}(k)-\hat{X}_{i}(k)\|_{p}^{2}+\tilde{\xi}_{11}(1+\varphi^{-1})\tilde U(k),
  \end{align}
where $\tilde{\xi}_{11}=5\hat{d}^{2}((8\gamma^{2}+2\eta^{2}L_{f}^{2})(1+L_{f}^{2})+2\eta L_{f}^{2}(2+\eta L_{f}^{2}))$.

Suppose that (\ref{Induction}) holds at $k$, we have
  \begin{align}\label{Inductionkx}
  &\|X_{i}(k)-\hat{X}_{i}(k)\|_{p}\notag\\
  &=\|X_{i}(k)-\hat{X}_{i}(k-1)-s(k)\mathcal{C}((X_{i}(k)-\hat{X}_{i}(k-1))/s(k))\|_{p}\notag\\
  &=s(k)\|(X_{i}(k)-\hat{X}_{i}(k-1))/s(k)-\mathcal{C}((X_{i}(k)-\hat{X}_{i}(k-1))/s(k))\|_{p}\notag\\
  &\leq(1-\varphi)s(k),
  \end{align}
where the first equality holds due to (\ref{Alg3a}) and (\ref{Alg3g}); and the inequality holds due to (\ref{Assum3}).

Then, from (\ref{ErrorxbarxnoE}), (\ref{xminushatx}), and (\ref{Inductionkx}), we have
\begin{align}\label{ErrorxbarxA3a}
&\|\bm{X}(k+1)-\bar{\bm{X}}(k+1)\|^{2}\notag\\
&\leq\delta\|\bm{X}(k)-\bar{\bm{X}}(k)\|^{2}+\frac{2\eta^{2}}{\gamma(1-\sigma)}\|\bm{Y}(k)-\bar{\bm{Y}}(k)\|^{2}+n\tilde{d}^2\frac{8\gamma}{1-\sigma}(1-\varphi)^{2}s^{2}(0)\mu^{2k}.
\end{align}

Similarly,
we have
  \begin{align}\label{Inductionky}
  &\|Y_{i}(k)-\hat{Y}_{i}(k)\|_{p}\notag\\
  &=\|Y_{i}(k)-\hat{Y}_{i}(k-1)-s(k)\mathcal{C}((Y_{i}(k)-\hat{Y}_{i}(k-1))/s(k))\|_{p}\notag\\
  &=s(k)\|(Y_{i}(k)-\hat{Y}_{i}(k-1))/s(k)-\mathcal{C}((Y_{i}(k)-\hat{Y}_{i}(k-1))/s(k))\|_{p}\notag\\
  &\leq(1-\varphi)s(k),
  \end{align}
where the first equality holds due to (\ref{Alg3b}) and (\ref{Alg3h}); and the inequality holds due to (\ref{Assum3}) and (\ref{Induction}).

Then, from (\ref{ErrorygnoE}), (\ref{xyminushaty}), and (\ref{Inductionky}), we have
\begin{align}\label{ErrorybaryA3a}
&\|\bm{Y}(k+1)-\bar{\bm{Y}}(k+1)\|^{2}\notag\\
&\leq(\delta+4\frac{2 L_{f}^{2}}{\gamma(1-\sigma)}\eta^{2})\|\bm{Y}(k)-\bar{\bm{Y}}(k)\|^{2}+4\frac{2 L_{f}^{2}}{\gamma(1-\sigma)}(4\gamma^{2}(1-\varphi)^{2}n\tilde{d}^2s^{2}(0)\mu^{2k}\notag\\
&\quad+(8\gamma^{2}+2\eta^{2} L_{f}^{2})\|\bm{X}(k)-\bar{\bm{X}}(k)\|^{2}+\eta^{2}\|\bar{\bm{g}}^{0}(k)\|^{2})+\frac{8\gamma}{1-\sigma}(1-\varphi)^{2}n\tilde{d}^2s^{2}(0)\mu^{2k}.
\end{align}
%
%
%


From (\ref{OptimalF2}), (\ref{ErrorxbarxA3a}), and (\ref{ErrorybaryA3a}), we have
\begin{align}\label{Convergence4D}
\tilde{U}(k+1)
&\leq\tilde{\theta}_{5}\|\bm{X}(k)-\bar{\bm{X}}(k)\|^{2}+\breve\theta_{13}\|\bm{Y}(k)-\bar{\bm{Y}}(k)\|^{2}+\tilde{\phi}n(F(\bar{X}(k))-F^{\star})\notag\\
&\quad-\tilde{\theta}_{6}\|\bar{\bm{g}}^{0}(k)\|^{2}-\theta_{9}\tilde{\phi}\|\bar{\bm{g}}(k)\|^{2}+2n\tilde{d}^2\frac{8\gamma}{1-\sigma}(1-\varphi)^{2}s^{2}(k)(1+2L_{f}^{2}),
\end{align}
where $\tilde{\theta}_{5}=\delta+\phi\xi_{1}+\tilde{\phi}\frac{\eta L_{f}^{2}}{2}$,~$\tilde{\theta}_{6}=\tilde{\phi}\frac{\eta}{4}-\phi\varepsilon_{1}$.

From (\ref{omega1a}) and (\ref{breve_omega7}), we have
\begin{equation}\label{tildeomeg1}
0<\tilde{\theta}_{5}<1-0.59(1-\sigma)\gamma=1-\breve\theta_{3}<1.
\end{equation}


From $\eta<\frac{(1-\sigma)^{2}\gamma}{40L_{f}}$, we have
\begin{equation}\label{tildeomeg2}
2\nu\tilde{\theta}_{6}=2\nu(\frac{0.1\gamma(1-\sigma)}{L_{f}^{2}}-\frac{8L_{f}^{2}\eta^{2}}{(1-\sigma)\gamma}\phi)>(\frac{32\gamma}{1-\sigma}-\frac{8\gamma}{1-\sigma})2\nu\phi=\frac{48\nu\gamma\phi}{1-\sigma}>0.
\end{equation}

From $\gamma<\frac{2L_{f}^{2}}{(1-\sigma)\nu}$, we have
\begin{equation}\label{tildeomeg2}
2\nu\tilde{\theta}_{6}=2\nu(\frac{0.1\gamma(1-\sigma)}{L_{f}^{2}}-\frac{8L_{f}^{2}\eta^{2}}{(1-\sigma)\gamma}\phi)<\frac{0.2\nu\gamma(1-\sigma)}{L_{f}^{2}}<0.4.
\end{equation}

From (\ref{barg0k}), (\ref{breveomeg2}), (\ref{Convergence4D}), (\ref{tildeomeg1}), (\ref{tildeomeg2}), and $\mu\in[\sqrt{\tilde{\theta}_{1}},1)$, we have
\begin{align}\label{Convergence4E}
\tilde{U}(k+1)
&<(1-\tilde{\theta}_{3})\breve{U}(k)+ 2(1+2L_{f}^{2})\frac{8n\tilde{d}^2(1-\varphi)^{2}}{1-\sigma}\gamma s^{2}(k)\notag\\
&<((1-\tilde{\theta}_{3})\tilde{\xi}_{5}+\tilde{\theta}_{2}\gamma)s^{2}(k)=\frac{\tilde{\theta}_{1}}{\mu^{2}}\tilde{\xi}_{5}s^{2}(k+1)\leq\tilde{\xi}_{5}s^{2}(k+1).
\end{align}

From (\ref{tildeomeg2}), we have
\begin{equation}\label{breveomega7}
1>\tilde{\theta}_{3}=\min\{\breve\theta_{3},~\frac{48\nu\gamma\phi}{1-\sigma}\}=\tilde{\theta}_{4}\gamma.
\end{equation}

From (\ref{breveomega7}), and $\tilde{\xi}_{5}>\frac{\tilde{\theta}_{2}}{\tilde{\theta}_{4}}$, we have
\begin{equation}\label{ineb}
0<\tilde{\theta}_{1}=1-\tilde{\theta}_{3}+\frac{\tilde{\theta}_{2}}{\tilde\xi_{5}}\gamma=1-(\tilde{\theta}_{4}-\frac{\tilde{\theta}_{2}}{\tilde{\xi}_{5}})\gamma<1.
\end{equation}

Based on (\ref{xhatxD3}), and (\ref{Inductionkx}), we have
  \begin{align}\label{xhatxD5}
  &\max_{i\in\mathcal{V}}\|X_{i}(k+1)-\hat{X}_{i}(k)\|_{p}^{2}\notag\\
  &\leq(1+\varphi+\tilde{\xi}_{8}\gamma^{2})(1-\varphi)^{2}s^{2}(k)+4\hat{d}^{2}(8\gamma^{2}+\eta^{2}(1+4L_{f}^{2})+4\eta)(1+\varphi^{-1})\tilde{\xi}_{5}s^{2}(k)\notag\\
  &<(1-(\varphi+\varphi^{2}-\varphi^{3})+\tilde{\xi}_{8}\gamma^{2}(1-\varphi)^{2}+4\hat{d}^{2}(8\gamma^{2}+\eta(5+4L_{f}^{2}))(1+\varphi^{-1})\tilde{\xi}_{5})s^{2}(k)\notag\\
  &=\frac{\tilde{\xi}_{6}}{\mu^{2}}s^{2}(k+1)\leq s^{2}(k+1),
  \end{align}
where the second inequality holds due to $\eta<1$; the last inequality holds due to $\eta<\frac{\varphi+\varphi^{2}-\varphi^{3}}{2\tilde{\xi}_{1}}$, $\gamma<\sqrt{\frac{\varphi+\varphi^{2}-\varphi^{3}}{2\tilde{\xi}_{3}}}$, and $\mu\in[\sqrt{\tilde{\xi}_{6}},1)$.


%
%
%

Similarly, we have
  \begin{align}\label{yhatyD5}
  &\max_{i\in\mathcal{V}}\|Y_{i}(k+1)-\hat{Y}_{i}(k)\|_{p}^{2}\notag\\
  &\leq((1+\varphi+\tilde{\xi}_{9}\gamma^{2}(1+L_{f}^{2}))(1-\varphi)^{2}+5\hat{d}^{2}((8\gamma^{2}+2\eta^{2}L_{f}^{2})(1+L_{f}^{2})\notag\\
  &\quad+2\eta L_{f}^{2}(2+\eta L_{f}^{2}))(1+\varphi^{-1})\tilde{\xi}_{5})s^{2}(k)\notag\\
  &<(1-(\varphi+\varphi^{2}-\varphi^{3})+\tilde{\xi}_{9}\gamma^{2}(1+L_{f}^{2})(1-\varphi)^{2}\notag\\
  &\quad+(40\hat{d}^{2}\gamma^{2}(1+L_{f}^{2})+10L_{f}^{2}\eta(3+2L_{f}^{2}))(1+\varphi^{-1})\tilde{\xi}_{5})s^{2}(k)\notag\\
  &=\frac{\tilde{\xi}_{7}}{\mu^{2}}s^{2}(k+1)\leq s^{2}(k+1),
  \end{align}
where the second inequality holds due to $\eta<1$; the last inequality holds due to $\eta<\frac{\varphi+\varphi^{2}-\varphi^{3}}{2\tilde{\xi}_{2}}$, $\gamma<\sqrt{\frac{\varphi+\varphi^{2}-\varphi^{3}}{2\tilde{\xi}_{4}}}$, and $\mu\in[\sqrt{\tilde{\xi}_{7}},~1)$.

Therefore, from (\ref{xhatxD5}) and (\ref{yhatyD5}), we know that (\ref{Induction}) holds at $k+1$. Finally, by mathematical induction, we know that (\ref{Induction}) holds for any $k>0$. Hence, (\ref{thm-7a}) holds.
%
%
This completes the proof.\\

\end{document}